\documentclass[12pt]{amsart} 
\usepackage{epsfig,amsthm,amsmath,amsfonts,amssymb}
\makeatletter 

\usepackage{graphics,graphicx,color} 
\usepackage{amssymb, amsmath, amsbsy,mathtools}

\newtheorem{theorem}{Theorem}[section]
\newtheorem{lemma}[theorem]{Lemma}
\newtheorem{corollary}[theorem]{Corollary}

\newtheorem{remark}[theorem]{Remark}

\def\proof{\noindent{\bf Proof.} }
\def\endproof{{}\hfill{$\square$}\par\medskip}

\newcommand{\ei}{\boldsymbol{e}_i}
\newcommand{\ej}{\boldsymbol{e}_j}
\newcommand{\oneomega}{\boldsymbol{1}_{\omega}}
\newcommand{\sHoneper}{\sH^1_{\operatorname{per}}}

\newcommand{\dd}{\mathrm{d}}
\newcommand{\DD}{\mathrm{D}}

\newcommand{\R}{\mathbb{R}}

\newcommand{\h}{\boldsymbol{h}}

\newcommand{\n}{\boldsymbol{n}}

\newcommand{\nablat}{{{\nabla}_{\boldsymbol\tau}}}

\newcommand{\sH}{{\rm H}}

\newcommand{\sL}{{\rm L}}

\renewcommand{\div}{{\operatorname{div}}}
\newcommand{\divt}{{\operatorname{div}_{\boldsymbol{\tau}}}}

\begin{document}
\title[Shape optimization for composite materials
and scaffolds]{Shape optimization for composite materials
and scaffold structures}
\author{Marc Dambrine}
\address{Marc Dambrine, 
CNRS / Univ Pau \& Pays Adour,  Laboratoire de Math\'ematiques et 
de leurs Applications de Pau -- F\'ed\'eration IPRA, UMR 5142, 64000, Pau, France.}
\email{marc.dambrine@univ-pau.fr}
\author{Helmut Harbrecht}
\address{Helmut Harbrecht,
Departement f\"ur Mathematik und Informatik, 
Universit\"at Basel, 
Spiegelgasse 1, 4051 Basel, Switzerland.}
\email{helmut.harbrecht@unibas.ch}

\begin{abstract}
This article combines shape optimization and homogenization
techniques by looking for the optimal design of the microstructure
in composite materials and of scaffolds. The development of 
materials with specific properties is of huge practical interest, for 
example, for medical applications or for the development of light weight 
structures in aeronautics. In particular, the optimal design of microstructures 
leads to fundamental questions for porous media: what is the sensitivity of 
homogenized coefficients with respect to the shape of the microstructure? 
We compute Hadamard's shape gradient for the problem of realizing 
a prescribed effective tensor and demonstrate the applicability and 
feasibility of our approach by numerical experiments.
\end{abstract}

\maketitle

\section{Introduction}
Shape optimization has been developed as an efficient
method for designing devices, which are optimized with 
respect to a given purpose. Many practical problems from engineering 
amount to boundary value problems for an unknown function, which 
needs to be computed to obtain a real quantity of interest. For example, 
in structural mechanics, the equations of linear elasticity are usually 
considered and solved to compute, for example, the leading mode of 
a structure or its compliance. Shape optimization is then applied to 
optimize the workpiece under consideration with respect to this objective 
functional, see \cite{ZOL1,MS,SI,PIR,ZOL2} and the references therein for 
an overview on the topic of shape optimization, which falls into the 
general setting of optimal control of partial differential equations.

In the present article, we will consider a slightly different question: 
the optimal design of microstructures in composite materials. Indeed, 
the additive manufacturing allows to build lattices or perforated structures 
and hence to build structures with physical properties that vary in space.
The realization of composite materials or -- as a limit case -- scaffold 
structures with specific properties has, of course, a huge impact for 
many practical applications. Examples arise from the development 
of light weight structures in aeronautics or for medical implants 
in the orthopedic and dental fields, see e.g.\ \cite{NRD,WXZ} 
and the references therein.

The optimal design of composite materials and 
scaffold structures has been considered in many works, 
see \cite{AOTHH,CCFK,FCHO,HMT,HWGDFS,LKH,LRC,
PFM,RSZ,S,SSW,WK,WWSK} for some of the respective 
results. The methodology used there is primarily based on 
the forward simulation of the material properties of a given 
microstructure. Whereas, sensitivity analysis has been 
used in \cite{AGP,GAP} to compute the derivatives with 
respect to the side lengths and the orientation of a rectangular 
inclusion. In \cite{HRLS}, the derivatives with respect to the coefficients 
of a B-spline parametrization of the inclusion have been computed.
In \cite{NC}, the shape derivative has been derived in the context 
of a level set representation of the inclusion. We are, however, not 
aware on optimization results which employ Hadamard's shape 
gradient \cite{HADA}. Therefore, in the present article, we perform
the sensitivity analysis of the effective material properties with 
respect to the shape of the inclusions: we compute the 
related shape gradient and consider its efficient computation 
by {\em homogenization}. As an application of these computations, 
we focus on the least square matching of a desired material property.

This article is structured as follows. In Section~\ref{sec:problem},
we briefly recall the fundamentals of homogenization theory and
introduce the problem under consideration. Then, Section
\ref{sec:calculus} is dedicated to shape calculus for composite
materials that are the mixture of two materials with different 
physical properties. We compute the local shape derivative for the cell 
functions and study the sensitivity of the effective tensor with 
respect to the microstructure. These results are extended to
scaffold structures in Section~\ref{sec:scaffolds}. Here, we also
provide second order shape derivatives, which can especially
be used in uncertainty quantification. Finally, we present 
numerical results in Section~\ref{sec:results} for the least 
square matching of the effective tensor. In particular, we 
exhibit various solutions for the same tensor in order to 
illustrating the non-uniqueness of the shape optimization 
problem under consideration.
 
\section{The problem and the notation}\label{sec:problem}
\subsection{Homogenization}
To describe the goals and methods of the present article, we 
shall restrict ourselves to the situation of the simple two-scale 
problem posed in a domain $D\subset\mathbb{R}^d$, $d=2,3$:
\begin{equation}\label{homo1}
-\div\big(\boldsymbol{A}^\varepsilon\nabla u^\varepsilon\big) = f\ \text{in $D$},
\quad u^\varepsilon=0\ \text{on $\partial D$}.
\end{equation}
Here, the $(d\times d)$-matrix $\boldsymbol{A}^\varepsilon$ is 
assumed to be oscillatory in the sense of 
\[
\boldsymbol{A}^\varepsilon(x) = \boldsymbol{A}\biggl({x\over\varepsilon}\biggr),
	\quad x\in D.
\]
Mathematical homogenization is the study of the limit of 
$u^\varepsilon$ when $\varepsilon$ tends to $0$. Various 
approaches have been developed to this end. The oldest one 
is comprehensively exposed in Bensoussan, Lions and Papanicolaou 
\cite{BLP}. It consists in performing a formal multiscale asymptotic expansion 
and then in the justification of its convergence using the energy 
method due to Tartar \cite{T}. A significant result obtained 
with this approach was the existence of a ($L^2(D)$-) limit 
$u_0(x)$ of $u^{\varepsilon}(x)$ and, more importantly, the 
{\em identification of a limiting, ``effective'' or ``homogenized'' 
elliptic problem in $D$ satisfied by $u_0$}.

\begin{figure}[hbt]
\begin{center}
\setlength{\unitlength}{0.8cm} 
\begin{picture}(13,7.2)
\put(0.015,0.02){\includegraphics[width=4.96cm,height=4.96cm]{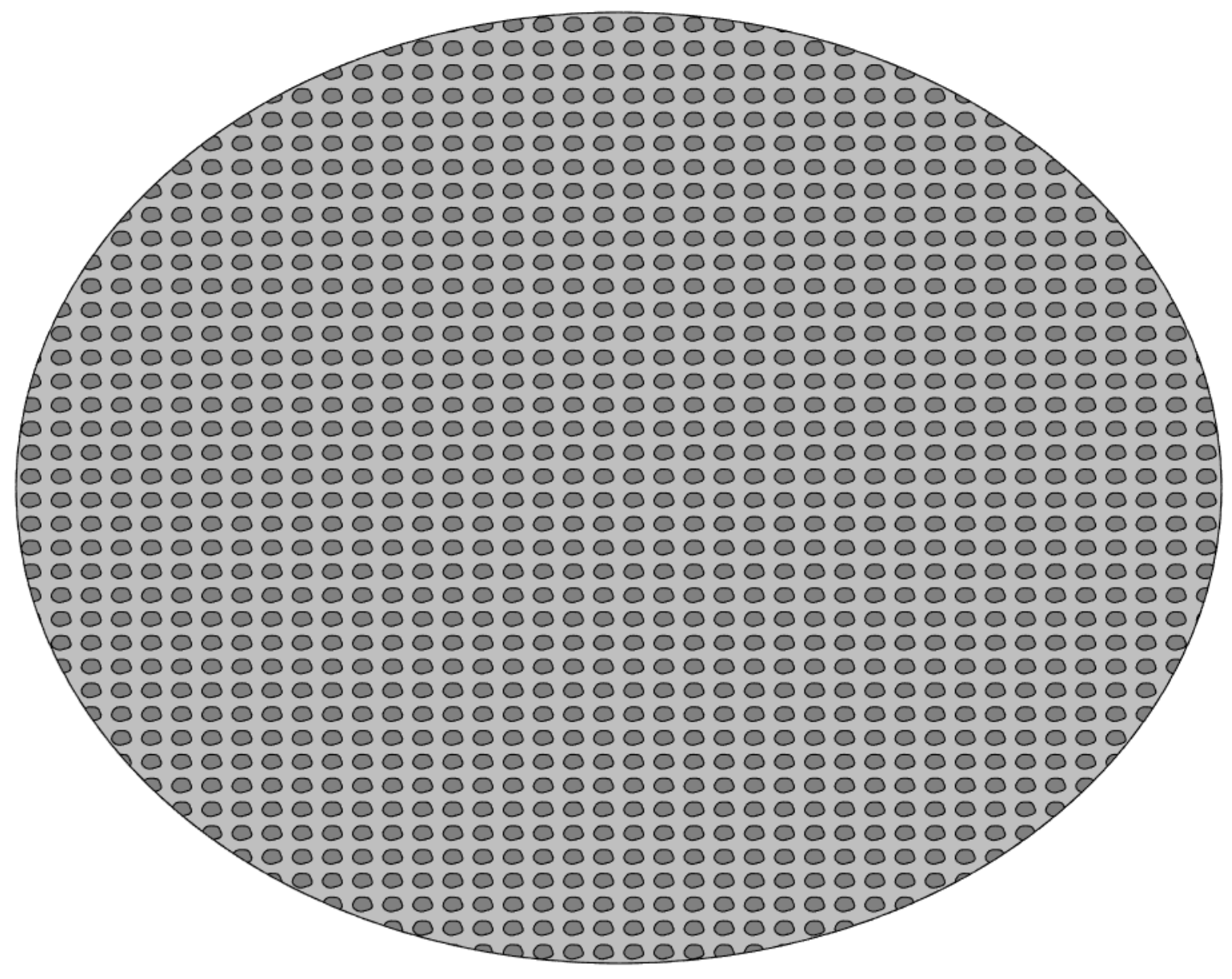}}
\put(2.975,2.975){\line(1,0){0.15}}
\put(2.975,2.975){\line(0,1){0.15}}
\put(2.975,3.125){\line(1,0){0.15}}
\put(3.125,2.975){\line(0,1){0.15}}
\put(2.975,3.125){\line(3,2){4.23}}
\put(2.975,2.975){\line(3,-2){4.23}}
\put(7.2,0.13){\includegraphics[width=4.66cm,height=4.66cm]{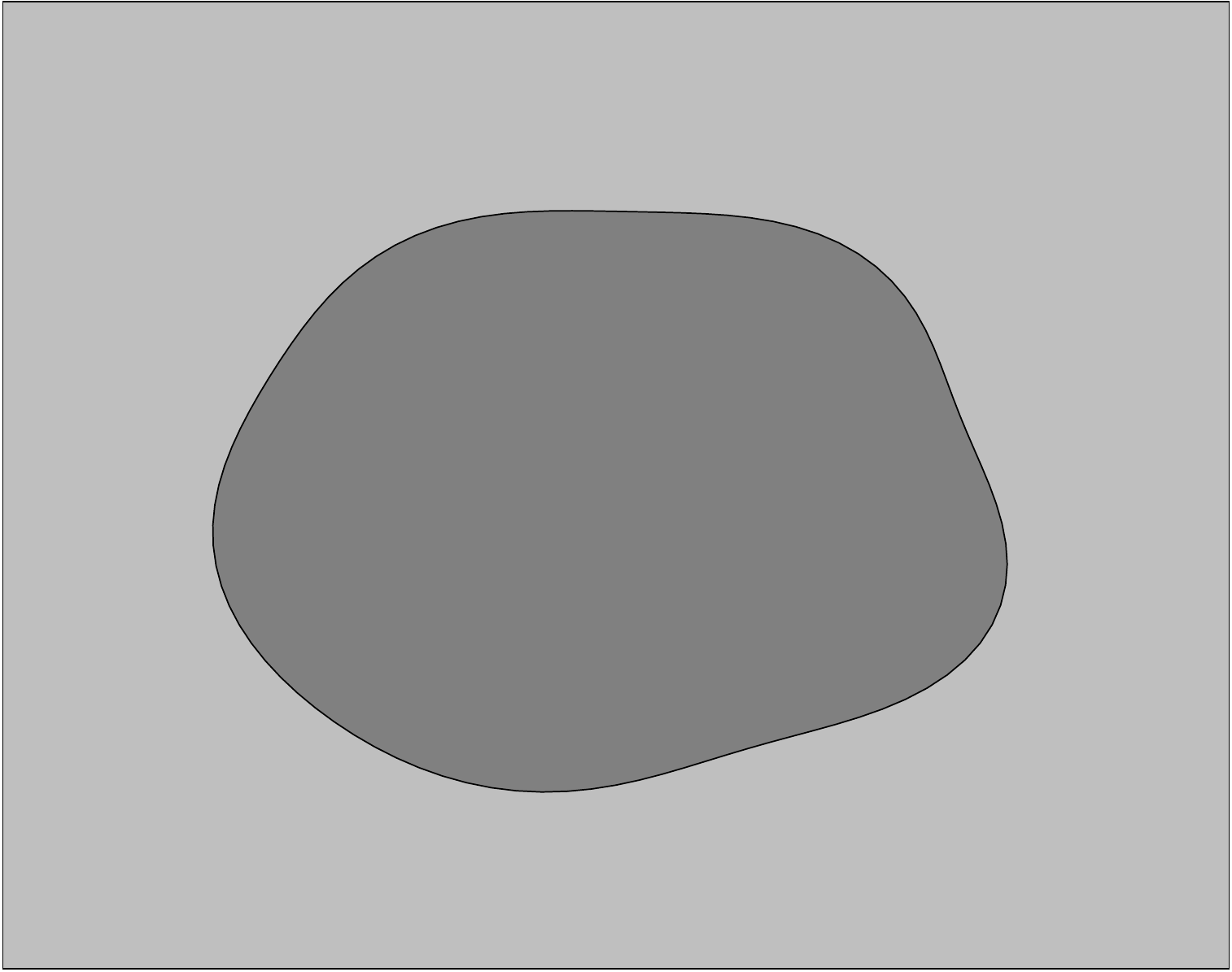}}
\qbezier[1000](3.125,3.125)(7.2,4.3)(7.2,4.3)
\qbezier[1000](3.125,2.975)(7.2,1.8)(7.2,1.8)
\put(7.2,0.13){\line(0,1){5.8}}
\put(7.2,0.13){\line(1,0){5.8}}
\put(13,0.13){\line(0,1){5.8}}
\put(7.2,5.95){\line(1,0){5.8}}
\put(0.5,5.2){$D$}
\put(7.3,5.5){$\partial Y$}
\put(9,4){$\partial\omega$}
\put(11.9,4.9){$Y$}
\end{picture}
\caption{
The domain $D$ with unit cell ${Y}$.}
\label{fig:setup}
\end{center}
\end{figure}
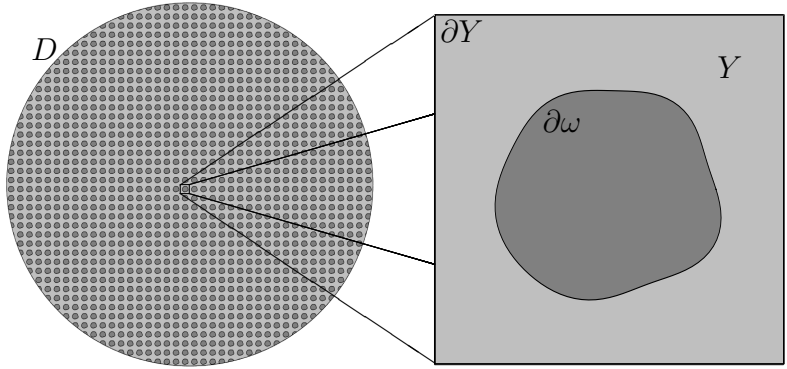

We introduce the unit cell $Y=[0,1]^d$ for the fast scale 
of the problem \eqref{homo1} and assume that the matrix 
function $\boldsymbol{A}(y)$ has period $Y$, cf.~Figure
\ref{fig:setup} for a graphical illustration. Moreover, we 
consider the space $\sHoneper(Y)$ of $Y$-periodic 
functions that belong to $\sH^1(Y)$ and the unit vector 
$\ei \in \R^d$ in the $i$-th direction of $\R^d$. Then, we 
can define the cell problems for all $i=1,\dots,d$:
\[
\text{find $w_i\in\sHoneper(Y)$ such that}
-\div\big(\boldsymbol{A}(\ei+\nabla w_i)\big) = 0.
\]
The Lax-Milgram theorem ensures the existence and 
uniqueness of the solutions $w_i$ to these cell problems
for $1\leq i \leq d$.

The family of functions $w_i$ can be used to define the 
\emph{effective tensor} $\boldsymbol{A}_0 = [a_{i,j}]_{i,j=1}^d$ 
in accordance with
\[
a_{i,j} = \int_Y \langle\boldsymbol{A}(\ei +\nabla w_i),\ej + \nabla w_j\rangle\dd y.
\]
It yields the \emph{homogenized solution\/} $u_0\in
\sH_0^1(D)$ by means of the limiting equation
\[
-\div\big(\boldsymbol{A}_0\nabla u_0\big)  = f\ \text{in}\ D,
\quad u_0 = 0\ \text{on}\ \partial D.
\]
In particular, by setting
\[
  u_1(x,y) = \sum_{i=1}^d \frac{\partial u_0}{\partial x_i}(x)w_i(y),
  	\quad (x,y)\in D\times Y,
\]
one has the error estimate
\[
  \bigg\|u^\varepsilon(x)-u_0(x)-\varepsilon u_1\bigg(x,{x\over\varepsilon}\bigg)\bigg\|_{H^1(D)}
  	\lesssim \sqrt{\varepsilon}\to 0\ \text{as $\varepsilon\to 0$},
\]
cf.~\cite{A,N89}.

\subsection{Composite materials}
From now on, we shall consider a composite material which
consists of two materials having a different conductivity.
Let $\omega$ be open subset of $Y$ and let 
\begin{equation}\label{eq:conductivity}
\sigma(y) = \sigma_1(y) + \big(\sigma_2(y)-\sigma_1(y)\big) \oneomega(y)
\end{equation}
be a piecewise smooth function defined on $Y$, where
$\sigma_1$ and $\sigma_2$ are smooth functions on $Y$
such that there exist real numbers $\overline{\sigma}
>\underline{\sigma}>0$ satisfying
\begin{equation}\label{eq:bounds}
\underline{\sigma}\leq \sigma_i(y)\leq\overline{\sigma}
\ \text{for all}\ y\in Y.
\end{equation}
The current situation is illustrated in Figure~\ref{fig:setup}.

In the following, we orient the surface $\partial \omega$ of the
inclusion $\omega$ so that its normal vector $\n$ indicates the 
direction going from the interior of $\omega$ to the exterior 
$Y\setminus\overline{\omega}$. The jump of a quantity $f$ 
through the interface $\partial\omega$ at a point 
$y\in\partial\omega$ is then
$$
[f(y)] = \lim_{t\rightarrow 0^+} f\big(y+t\n(y)\big)
	- \lim_{t\rightarrow 0^-} f\big(y+t\n(y)\big).
$$ 

For our model problem, the entries of the effective tensor 
$\boldsymbol{A}_0$ are given by
\begin{equation}\label{coef:effective:tensor}
a_{i,j}(\omega) =\int_Y \sigma\langle\ei +\nabla w_i,\ej + \nabla w_j\rangle\dd y,
\end{equation}
where the $w_i$ solves the respective cell problem
\begin{equation}\label{cell:problem}
\text{find $w_i\in\sHoneper(Y)$ such that}\ 
	-\div\big(\sigma(\ei+\nabla w_i)\big) = 0.
\end{equation}
Notice that $\boldsymbol{A}_0$ is a symmetric matrix, but it is in 
general not the identity, since the geometric inclusion generates 
anisotropy.

In this article, we consider a given tensor $\boldsymbol{B}
\in\mathbb{R}_{\text{sym}}^{d\times d}$ describing the desired 
material properties. We then address the following question: 
\emph{can we find a mixture (that is a domain $\omega$) 
such that the effective tensor is as close as possible 
to $\boldsymbol{B}$?} 
 
In order to make precise the notion of closeness 
between matrices, we choose the Frobenius norm 
on matrices and define the objective $J(\omega)$ 
to minimize
\begin{equation}\label{eq:functional}
J(\omega) = \cfrac{1}{2}\|\boldsymbol{A}_0(\omega)-\boldsymbol{B}\|_{F}^2 
	= \frac{1}{2} \sum_{1\leq i,j\leq d} \big(a_{i,j}(\omega)-b_{i,j}\big)^2.
\end{equation}

Of course, not every tensor can be reached by mixing 
two materials. There exist bounds for the eigenvalues of 
the tensor $\boldsymbol{A}_0$. For example, the Voigt--Reuss 
bounds state that
$$
  \left(\int_Y \sigma^{-1}\dd y \right)^{-1}\|\boldsymbol{t}\|^2
	\leq \langle\boldsymbol{A}_0(\omega)\boldsymbol{t},\boldsymbol{t}\rangle
	\leq\left(\int_Y \sigma\dd y\right)\|\boldsymbol{t}\|^2,
$$
compare \cite{Reuss,Voigt}. Hence, the infimum is not $0$ 
if the target tensor is not in the closure of tensors reachable 
by a mixture of two materials. 

\section{Shape calculus for the mixture}\label{sec:calculus}
\subsection{Local shape derivative}
We introduce a vector field $\h:Y\to Y$ that vanishes on the 
boundary $\partial Y$ of the reference cell but whose 
action may deform the interior interface $\partial\omega$. 
We consider the perturbation of identity $T_t=I+t\h$, which
is a diffeomorphism for $t$ small enough that preserves 
$Y$. We denote by $\omega(t) =T_t(\omega)$, $\sigma(t,y) 
= \sigma\big(T_t(y)\big)$, and $w_i(t)\in\sHoneper(Y)$ the 
solution of \eqref{cell:problem} for the inclusion $\omega(t)$. 

We are interested in describing how the effective tensor 
depends on the deformation field $\h$. We will successively 
consider the sensitivity on $\h$ first of the solutions $w_i(t)$ to the
cell problems, then of each entry of the effective conductivity tensor, 
and finally of the least square matching to the desired tensor. As it
turns out, it is much more convenient to compute the sensitivity
with respect to $w_i(t)$ in an indirect way by considering the 
shape derivative of the function $\phi_i = w_i+x_i$.

\begin{lemma}[Shape derivative of the cell problem]
\label{lemma:shape:derivative:cell:function}
The shape derivative $\phi_i'\in\sHoneper(Y)$ of the function 
$\phi_i = w_i+x_i$ is the solution in $Y$ of the transmission 
problem
\begin{equation}\label{shape:derivative:cell:problem}
\begin{aligned}
\Delta \phi_i'&= 0 &&\text{ in }\omega \cup Y\setminus\overline{\omega},\\[0.5ex]
[\phi'_i] &= -\langle\h,\n\rangle [\partial_{\n} \phi_i] &&\text{ on }\partial\omega,\\[0.5ex]
[\sigma\partial_{\n} \phi'_i]&= [\sigma]\,\divt(\langle\h,\n\rangle\nablat \phi_i) 
	&&\text{ on }\partial \omega.
\end{aligned}
\end{equation}
\end{lemma}

\proof
We proceed in the usual elementary way. Prove existence of the 
material derivative, then characterize it, and finally express the 
local shape derivative. 

\textsl{First step: computing the material derivative.} 
Let us prove that the material derivative exists and satisfies 
the variation problem: find $\dot{\phi}_i\in\sHoneper(Y)$ such that
\begin{equation}\label{weak:form:material:derivative}
\int_Y \sigma\langle\nabla\dot{\phi}_i,\nabla v\rangle\dd y 
	= - \int_Y \sigma\langle\big(\div \h I -(\DD\h^t+\DD\h)\big)\nabla\phi_i,\nabla v\rangle\dd y
\end{equation}
for all $v \in \sHoneper(Y)$. The transported function 
$\widetilde{\phi}_i(t,y)=\phi_i\big(t,T_t(y)\big)$ satisfies the variational equation
\begin{equation}\label{eq:transported}
\int_Y \sigma(y) \langle A(t,y)\nabla\big(\widetilde{\phi}_i(t,y)+\ei\big),\nabla v(y)\rangle\dd y = 0
\end{equation}
for all $v \in \sHoneper(Y)$, where we have set
\[
  A(t,y) = \DD T_t^{-1}(y)\big(\DD T_t^{-1}(y)\big)^t\det\big(\DD T_t(y)\big).
\]
We subtract from equation \eqref{eq:transported} the equation 
satisfied by $\phi_i$ for the reference configuration
\[
  \int_Y \sigma(y)\langle\nabla\big(\phi_i(y)+\ei\big),\nabla v(y)\rangle\dd y = 0
\]
to get for any $t>0$
\begin{equation}\label{diff:weak:form}
\begin{aligned}
  &\int_Y \bigg\{\sigma(y) A(t,y)\Big\langle\cfrac{\nabla \widetilde{\phi}_i(t,y)-\nabla \phi_i(y)}{t},\nabla v(y)\Big\rangle\\
  &\hspace*{3cm}+ \sigma(y)\Big\langle\cfrac{A(t,y)-I}{t}\nabla \phi_i(y),\nabla v\Big\rangle\bigg\}\dd y = 0
\end{aligned}
\end{equation}
for all $v \in \sHoneper(Y)$. Using $\widetilde{\phi}_i(t,\cdot)-\phi_i$ 
as test function and observing \eqref{eq:bounds}, we obtain 
the upper bound
\[
  \frac{\underline{\sigma}}{2}\left\| \cfrac{\nabla\widetilde{\phi}_i(t,\cdot)-\nabla\phi_i}{t}\right\|_{\sL^2(Y)}
  \leq \left\|\cfrac{A(t,\cdot)-I}{t}\right\|_{\sL^{\infty}(Y)} \|\nabla\phi_i\|_{\sL^2(Y)}.
\]
Since the fraction $\big(\widetilde{\phi}_i(t,\cdot)-\phi_i\big)/t$ is bounded 
in $\sHoneper(Y)$, it is weakly convergent in $\sHoneper(Y)$ and 
its weak limit is the material derivative $\dot{\phi}_i \in\sHoneper(Y)$. 

In order to prove the strong convergence of  $\big(\widetilde{\phi}_i(t,\cdot)-\phi_i\big)/t$ 
to  $\dot{\phi}_i \in\sHoneper(Y)$, we use $\big(\widetilde{\phi}_i(t,\cdot)-\phi_i\big)/t$ as 
test function in \eqref{diff:weak:form}:
\begin{align*}
  &\int_Y \sigma(y) \Big\langle A(t,y)\cfrac{\nabla\widetilde{\phi}_i(t,x)-\nabla\phi_i(x)}{t},
  	\cfrac{\nabla\widetilde{\phi}_i(t,y)-\nabla\phi_i(y)}{t}\Big\rangle\dd y\\
  &\hspace*{2cm} = -\int_Y \sigma(y)\Big\langle\cfrac{A(t,y)-I}{t}\nabla\phi_i(y),
  	\cfrac{\nabla\widetilde{\phi}_i(t,y)-\nabla\phi_i(y)}{t}\Big\rangle\dd y.
\end{align*}
We split the right-hand side into $-\big(R_1(t)+R_2(t)\big)$,
where
$$
R_1(t)= \int_Y \sigma(y)\Big\langle \big(A(t,y)-I\big)
	\cfrac{\nabla\widetilde{\phi}_i(t,y)-\nabla\phi_i(y)}{t},
		\cfrac{\nabla\widetilde{\phi}_i(t,y)-\nabla\phi_i(y)}{t}\Big\rangle\dd y
$$
and
$$		
R_2(t)= \int_Y \sigma(y)\Big\langle\cfrac{A(t,y)-I}{t}\nabla \widetilde{\phi}_i(t,y),
	\cfrac{\nabla\widetilde{\phi}_i(t,y)-\nabla\phi_i(y)}{t}\Big\rangle\dd y.
$$
The weak convergence of $\big(\widetilde{\phi}_i(t,\cdot)-\phi_i\big)/t$ 
to  $\dot{\phi}_i\in\sHoneper(Y)$ amounts to $R_1(t)\to 0$ 
while
\[
  R_2(t)\to\int_Y\sigma(y)\Big\langle
  	\left(\cfrac{\dd A(t,y)}{\dd t}\right)_{|t=0}\nabla\phi_i(y),\nabla\dot{\phi}_i(y)\Big\rangle\dd y.
\]  
Let us recall that it follows by an elementary computation that
\[
\left(\cfrac{\dd A(t,\cdot)}{\dd t}\right)_{|t=0} 
	= \div{\h} I -\left(\DD\h^t+\DD\h \right):=\mathcal{A}.
\]
We hence conclude that $\nabla\big(\widetilde{\phi}_i(t,\cdot)-\phi_i\big)/t$ 
converges strongly in $\sL^2(Y)$, which implies by Poincar\'e's 
inequality that $\big(\widetilde{\phi}_i(t,\cdot)-\phi_i\big)/t$ converges 
strongly in $\sHoneper(Y)$ to the material derivative.

\textsl{Second step: recovering the shape derivative.} 
We recall that the shape derivative is obtained from the material 
derivative by the relationship $\phi'_i=\dot{\phi}_i - \langle\h,\nabla\phi_i\rangle$. 
The first transmission condition $[\phi'_i] = -\langle\h,\n\rangle 
[\partial_{\boldsymbol{n}}\phi_i]$ on $\partial\omega$ expresses
that $\dot{\phi}_i$ has no jump on $\partial\omega$ as a function 
in $\sHoneper(Y)$.  

In order to prove the remaining relations, we decompose 
$\mathcal{A}$ to arrive at
\begin{align*}
\int_Y \sigma\langle\nabla \dot{\phi}_i,\nabla v\rangle\dd y 
	&= - \int_Y \sigma\langle\mathcal{A} \nabla\phi_i,\nabla v\rangle\dd y\\
	&= - \int_Y \sigma\big\{\div\h\,\langle\nabla\phi_i,\nabla v \rangle
	+\langle\h,\nabla\phi_i\rangle\Delta v +\langle\h,\nabla v\rangle\Delta\phi_i\big\}\dd y
\end{align*}
for all $v\in\sHoneper(Y)$. Next, we integrate by parts in 
$\omega$ and $Y\setminus\overline{\omega}$:
\begin{align*}
\int_Y \sigma\langle\nabla\dot{\phi}_i ,\nabla v\rangle\dd y
	&= \int_{\partial\omega} \big\{\left[\sigma\langle\nabla\phi_i,\nabla v\rangle\langle\h,\n\rangle\right] 
		- \left[\sigma\langle\h,\nabla v\rangle\partial_{\boldsymbol{n}}\phi_i\right]\big\}\dd o\\
	&\qquad+ \int_Y \sigma \big\langle\nabla\langle\h,\nabla\phi_i\rangle,\nabla v\big\rangle\dd y.
\end{align*}
This leads to
\begin{align*}
\int_Y \sigma\big\langle\nabla(\dot{\phi}_i - \langle\h,\nabla\phi_i\rangle),\nabla v\big\rangle\dd y
	&= \int_{\partial\omega}(\sigma_1-\sigma_2) \langle\h,\n\rangle\langle\nablat\phi_i,\nablat v\rangle\dd o \\
	&= \int_{\partial\omega}(\sigma_2-\sigma_1) \divt(\langle\h,\n\rangle\nablat\phi_i)v\dd o.
\end{align*}
We deduce that $\phi_i'$ is harmonic in both subdomains, $\omega$ 
and $Y\setminus\overline{\omega}$. Moreover, we obtain the second 
transmission condition.
\endproof

\subsection{Sensitivity of the effective tensor}
With the help of the local shape derivative \eqref{shape:derivative:cell:problem},
we can now compute the shape derivative of the effective tensor 
by using the basic formula of Hadamard's shape calculus
\begin{equation}\label{eq:rule}
\bigg(\cfrac{\dd}{\dd t} \int_{T_t(Y)} f(t,y)\,\dd y\bigg)\bigg|_{t=0} 
	=  \int_Y \bigg\{\cfrac{\dd}{\dd t} f(t,y) +\div\big(f(t,y)\h\big)\bigg\}\dd y,
\end{equation}
compare \cite{ZOL1,ZOL2} for example.

\begin{lemma} [Shape derivative of the coefficients of the effective tensor]
\label{lemma:shape:derivative:effective:tensor}
The shape derivatives of the entries $a_{i,j}$ of the effective tensor 
given by \eqref{coef:effective:tensor} are 
\begin{equation}\label{shape:derivative:entry:effective:tensor}
a_{i,j}'(\omega)[\h]= \int_{\partial\omega}[\sigma]\big\{
 \partial_{\n} \phi_i^- \partial_{\n} \phi_j^+ 
 +\langle\nablat \phi_i,\nablat \phi_j\rangle\big\}\langle\h,\n\rangle\,\dd o,
\end{equation}
where $\phi_i := w_i+x_i$ with $w_i\in \sHoneper(Y)$ being 
the solution to the $i$-th cell problem \eqref{cell:problem}.
\end{lemma}

\proof
We find
\begin{align*}
a_{i,j}'(\omega)[\h]&= \int_Y \sigma\big\{\langle\nabla \phi_i,\nabla \phi_j'\rangle
	+ \langle\nabla \phi'_i,\nabla \phi_j\rangle\big\}  \dd y 
	+ \int_Y \div(\sigma\langle\nabla \phi_i,\nabla \phi_j\rangle \h)\dd y\\
  &= \int_{\partial\omega} \big\{\phi_i [\sigma \partial_{\n} \phi_j'] +\phi_j [\sigma \partial_{\n} \phi_i']\big\}\dd o 
  	+ \int_{\partial\omega} [\sigma\langle\nabla \phi_i,\nabla \phi_j\rangle]\langle\h,\n\rangle\dd o.
\end{align*}
Moreover, one computes 
\begin{align*}
 [\sigma\langle\nabla\phi_i,\nabla\phi_j\rangle] 
 	&= \big[\sigma\big\{\langle\nablat\phi_i,\nablat\phi_j\rangle+\partial_{\n}\phi_i \partial_{\n}\phi_j\big\}\big]\\
 	&= [\sigma]\big\{\langle\nablat\phi_i,\nablat\phi_j\rangle-\partial_{\n}\phi_i^-\partial_{\n}\phi_j^+\big\}
\end{align*}
by employing the formula 
\[
  [abc] = [ab]c_1-[a]b_2c_1+[ac] b_2.
\] 
Hence, in view of the jump conditions $[\sigma \partial_{\n} \phi_i'] = 
[\sigma]\,\divt(\langle\h,\n\rangle\nablat\phi_i)$ and using integration 
by parts, we arrive at 
\begin{align*}
  a_{i,j}'(\omega)[\h]&=\int_{\partial\omega} [\sigma]\big\{
  	\phi_i\divt(\nablat\phi_j\langle\h,\n\rangle)+\phi_j\divt(\nablat\phi_i\langle\h,\n\rangle)\big\}\dd o\\
  	&\hspace*{2cm}+ \int_{\partial\omega}[\sigma]\big\{\langle\nablat\phi_i,\nablat\phi_j\rangle
	- \partial_{\n}\phi_i^-\partial_{\n}\phi_j^+\big\}\langle\h,\n\rangle\dd o\\
	&= -\int_{\partial\omega}[\sigma]\big\{\langle\nablat\phi_i,\nablat\phi_j\rangle
	+\partial_{\n}\phi_i^-\partial_{\n}\phi_j^+\big\}\langle\h,\n\rangle\Big\}\dd o.
\end{align*}
\endproof

\begin{remark}
By substituting back $\phi_i = w_i+x_i$, we
immediately arrive at the computable expression
\begin{align*}
a_{i,j}'(\omega)[\h]&= -\int_{\partial\omega}[\sigma]\big\{\langle\nablat (w_i+x_i),\nablat (w_j+x_j)\rangle\\
	&\hspace*{2cm}+(\partial_{\n} w_i^-+n_i)(\partial_{\n} w_j^++n_j)\big\}\langle\h,\n\rangle\dd o.
\end{align*}
Herein, $n_i$ denotes the $i$-th component of the normal vector $\n$.
\end{remark}

\subsection{Shape gradient of the least square matching}
With the help of Lemma~\ref{lemma:shape:derivative:effective:tensor}
and the chain rule
\[
  J'(\omega)[\h] = \sum_{1\leq i,j\leq d}\big(a_{i,j}(\omega)-b_{i,j}\big)
  	a_{i,j}'(\omega)[\h],
\]
we can easily determine the shape derivative of the 
objective $J(\omega)$ given by \eqref{eq:functional}.

\begin{corollary}
The shape derivative of the objective $J(\omega)$ from 
\eqref{eq:functional} reads
\begin{equation}\label{shape:derivative:objective}
 \begin{aligned}
 J'(\omega)[\h] &= \sum_{1\leq i,j\leq d} \big(b_{i,j}-a_{i,j}(\omega)\big)\\
	&\hspace*{1.5cm}\times\int_{\partial\omega}[\sigma]\big\{\langle\nablat\phi_i,\nablat\phi_j\rangle
		+\partial_{\n}\phi_i^-\partial_{\n}\phi_j^+\big\}\langle\h,\n\rangle\dd o.
\end{aligned}
\end{equation}
\end{corollary}

\section{Perforated plates and scaffold structures}\label{sec:scaffolds}
\subsection{Mathematical formulation}
We shall next consider the situation that $\sigma_2\to0$ in 
\eqref{eq:conductivity}, i.e., the case of perforated plates, which 
appear in two spatial dimensions, and scaffold structures, which 
appear in three spatial dimensions. In other words, the unit cell $Y = 
[0,1]^d$ contains a hole $\omega$. The collection of interior boundaries 
being translates of $\varepsilon\partial\omega$ of the macroscopic 
domain $D^\varepsilon$ is denoted by $\partial D_{\text{int}}^\varepsilon$
while the rest of the boundary $\partial D^\varepsilon\setminus\partial 
D_{\text{int}}^\varepsilon$ is denoted by $\partial D_{\text{ext}}^\varepsilon$.
In accordance with \cite{CIO}, we shall consider the boundary 
value problem
\begin{equation}\label{homo2}
\begin{aligned}
-\div\big(\boldsymbol{A}^\varepsilon\nabla u^\varepsilon\big) &= f\ &&\text{in $D^\varepsilon$},\\
\langle\boldsymbol{A}^\varepsilon\nabla u^\varepsilon,\n\rangle &=0\ 
&&\text{on $\partial D_{\text{ext}}^\varepsilon$},\\
u^\varepsilon&=0\ &&\text{on $\partial D_{\text{int}}^\varepsilon$}.
\end{aligned}
\end{equation}

To derive the homogenized problem, one introduces the cell functions 
$w_i\in\sHoneper(Y\setminus\overline{\omega})$ which are now 
given by the Neumann boundary value problems
\[
\begin{aligned}
  -\div\big(\sigma(\ei+\nabla w_i)\big) &= 0 &&\quad \text{in $Y\setminus\overline{\omega}$},\\
   \partial_{\n}w_i &= -n_i &&\quad \text{on $\partial\omega$}.
\end{aligned}
\]
The homogenized equation becomes
\[
-\div\big(\boldsymbol{A}_0\nabla u_0\big)  = (1-|\omega|) f\ \text{in}\ D,
\quad u_0 = 0\ \text{on}\ \partial D.
\]
Here, the domain $D$ coincides with $D^\varepsilon$ 
except for the holes and the effective tensor $\boldsymbol{A}_0 =
[a_{i,j}]_{i,j=1}^d$ is now given by
\begin{equation}\label{equi:tensor:Neumann}
a_{i,j}(\omega) =\int_{Y\setminus\overline{\omega}}
	\sigma\langle\ei +\nabla w_i,\ej + \nabla w_j\rangle\dd y,
\end{equation}
compare \cite{CIO}.

\subsection{Computation of the shape gradient}
Of course, the expression \eqref{equi:tensor:Neumann} looks like 
\eqref{coef:effective:tensor}, i.e., the one obtained in the case 
of a mixture. Indeed, this is normal since the case of a perforated 
domain can be seen as the limit case of the mixture when the 
inner conductivity $\sigma_2$ tends to $0$. This physically 
natural idea, purely isulation can be approximated by a very poorly 
conducting layer, has given birth to the ersatz material method in 
structural optimization \cite{AllaireJouveToader} and error estimates 
are given in \cite{DambrineKateb}. Nevertheless, due to the double 
passing to the limit, one cannot pass directly to the limit $\sigma_2 \rightarrow 0$ 
in the expression of the shape gradient of the mixture case given in 
Lemma \ref{lemma:shape:derivative:effective:tensor} to derive the 
shape gradient of \eqref{equi:tensor:Neumann}.

Therefore, we make again the ansatz $\phi_i = w_i+x_i$ and observe 
that $\phi_i$ satisfies homogeneous Neumann boundary conditions.
The derivative with respect to the shape in case of homogeneous 
Neumann boundary conditions is well-known. The shape derivative 
$\phi_i'$ of the state $\phi_i$ reads
\begin{equation}	 \label{locder}
\begin{aligned}
  \Delta\phi_i' &= 0 &&\quad \text{in $Y\setminus\overline{\omega}$},\\
   \partial_{\n}\phi_i' &= \divt(\langle\h,\n\rangle\nablat \phi_i)
  		&&\quad \text{on $\partial\omega$},
\end{aligned}
\end{equation}
see for example \cite{ZOL2} for the details of its computation.

In view of
\[
a_{i,j}(\omega)= \int_{Y\setminus\overline{\omega}}\sigma\langle\nabla\phi_i,
			\nabla\phi_j\rangle\dd y
\]
and \eqref{eq:rule}, the shape derivative of the coefficients of the 
effective tensor reads
\begin{align*}
a_{i,j}'(\omega)[\h]&= \int_{Y\setminus\overline{\omega}}\sigma\big\{\langle\nabla \phi_i,\nabla \phi_j'\rangle
	+ \langle\nabla \phi'_i,\nabla \phi_j\rangle\big\}  \dd y 
	+ \int_{Y\setminus\overline{\omega}} \div(\sigma\langle\nabla \phi_i,\nabla \phi_j\rangle \h)\dd y\\
  &= \int_{\partial\omega} \sigma\big\{\phi_i\partial_{\n} \phi_j' +\phi_j\partial_{\n} \phi_i'\big\}\dd o 
  	+ \int_{\partial\omega} \sigma\langle\nabla \phi_i,\nabla \phi_j\rangle\langle\h,\n\rangle\dd o.
\end{align*}
Inserting the boundary conditions of $\phi_i'$ and
integrating by parts gives
\begin{align*}
a_{i,j}'(\omega)[\h]&= \int_{\partial\omega} \sigma\big\{-2\langle\nablat\phi_i,\nablat \phi_j\rangle
  	+ \langle\nabla \phi_i,\nabla \phi_j\rangle\big\}\langle\h,\n\rangle\dd o\\
&= \int_{\partial\omega} \sigma\big\{\partial_{\n}\phi_i\partial_{\n}\phi_j
	-\langle\nablat\phi_i,\nablat \phi_j\rangle\big\}\langle\h,\n\rangle\dd o.
\end{align*}
Consequently, since it holds $\partial_{\n}\phi_i = 0$ on $\partial\omega$, 
the shape derivative of the objective $J(\omega)$ from \eqref{eq:functional} reads
\[
 \begin{aligned}
 J'(\omega)[\h] &= \sum_{1\leq i,j\leq d} \big(b_{i,j}-a_{i,j}(\omega)\big)
	\int_{\partial\omega}\sigma
		\langle\nablat\phi_i,\nablat\phi_j\rangle\langle\h,\n\rangle\dd o
\end{aligned}
\]
in the case of perforated plates or scaffold structures.

\subsection{Second order shape sensitivity analysis.} 
Using the expression of the entries of the effective tensor with 
respect to $\phi_i$, one checks that the diagonal entries are 
nothing but the Dirichlet energy associated to the solution of a 
homogeneous problem. The second order shape sensitivity 
analysis for this problem has been performed in 
\cite[Section 3]{DSZ}. We quote the result:
\begin{align*}
  a_{i,i}''(\omega)[\h,\h] &= \frac{1}{2}\int_{\partial\omega}
  	\big\{2\langle\nabla\phi_i', \nabla \phi_i \rangle
  + \div(\h)\|\nabla\phi_i\|^2 
  	+ \langle\nabla\|\nabla \phi_i\|^2,\h\rangle\big\}\langle \h,\n \rangle\,\dd o.
\end{align*}
Going carefully through the proofs presented there,
one immediately concludes that the off-diagonal terms are 
given by
\begin{align*}
 a_{i,j}''(\omega)[\h,\h] &= \frac{1}{2}\int_{\partial\omega}
 	\big\{\langle\nabla\phi_i', \nabla \phi_j \rangle 
		+ \langle\nabla\phi_i,\nabla\phi_j' \rangle \\ 
  &\qquad+ \div(\h)\langle\nabla\phi_i,\nabla\phi_j\rangle 
	+ \langle\nabla\langle\nabla\phi_i,\nabla\phi_j\rangle,\h\rangle\big\}\langle \h,\n \rangle\,\dd o.
\end{align*}

Second order shape derivatives can be used to
quantify uncertainties in the geometric definition of the 
microstructure, compare \cite{DHP,HAR}. Such uncertainties 
are motivated by tolerances in the fabrication process, for 
example by additive manufacturing. Manufactured devices are 
close to a nominal geometry but differ of course from there 
mathematical definition. Hence, the perturbations can be
assumed to be small.

The idea of the uncertainty quantification of the effective
tensor is as follows. For the perturbed inclusion 
$\omega_\varepsilon$, described by
\[
 \omega_\varepsilon = \{y+\varepsilon\h(y): y\in\omega\}\Subset Y, 
\]
we have the \emph{shape Taylor expansion}
\begin{align*}
 a_{i,j}(\omega_\varepsilon) &= a_{i,j}(\omega)
 	+ \varepsilon \alpha_{i,j}'(\omega)[\h] 
   	+ \frac{\varepsilon^2}{2} a_{i,j}''(\omega)[\h,\h] + \mathcal{O}(\varepsilon^3).
\end{align*}
Hence, first and even second order perturbation techniques can
be applied to quantify uncertainty in the geometric definition of
$\omega$. This means that, under the assumption that the 
perturbation $\h$ is a bounded random field and $\varepsilon$ 
is small, the expectation and the variance of the coefficient 
$a_{i,j}(\omega)$ of the effective tensor can be computed up 
to third order accuracy in the random perturbation's amplitude 
$\varepsilon$. We refer the interested reader to \cite{DHP,HAR} 
for all the details.

\section{Numerical results for the microstructure design}\label{sec:results}
\subsection{Implementation}
Our implementation is for the two-dimensional setting,
i.e., $d=2$. Especially, we will assume that the sought 
domain $\omega$ is starlike with respect to the midpoint 
of the unit cell. Then, we can represent its boundary 
$\partial\omega$ by using polar coordinates in 
accordance with
\[
  \partial\omega = \bigg\{r(\phi)\begin{bmatrix}\cos\phi\\\sin\phi\end{bmatrix}:
  \phi\in[0,2\pi)\bigg\},
\]
where the radial function $r(\phi)$ is given by the finite
finite Fourier series 
\begin{equation}\label{eq:Fourier}
  r(\phi) = a_0 + \sum_{k=1}^N \{a_k\cos(k\phi) + a_{-k}\sin(k\phi)\}.
\end{equation}
This yields the $2N+1$ design parameters $\{a_{-N},a_{1-N},\ldots,a_N\}$.

\begin{figure}[hbt]
\begin{center}
\includegraphics[width=11cm,height=5cm]{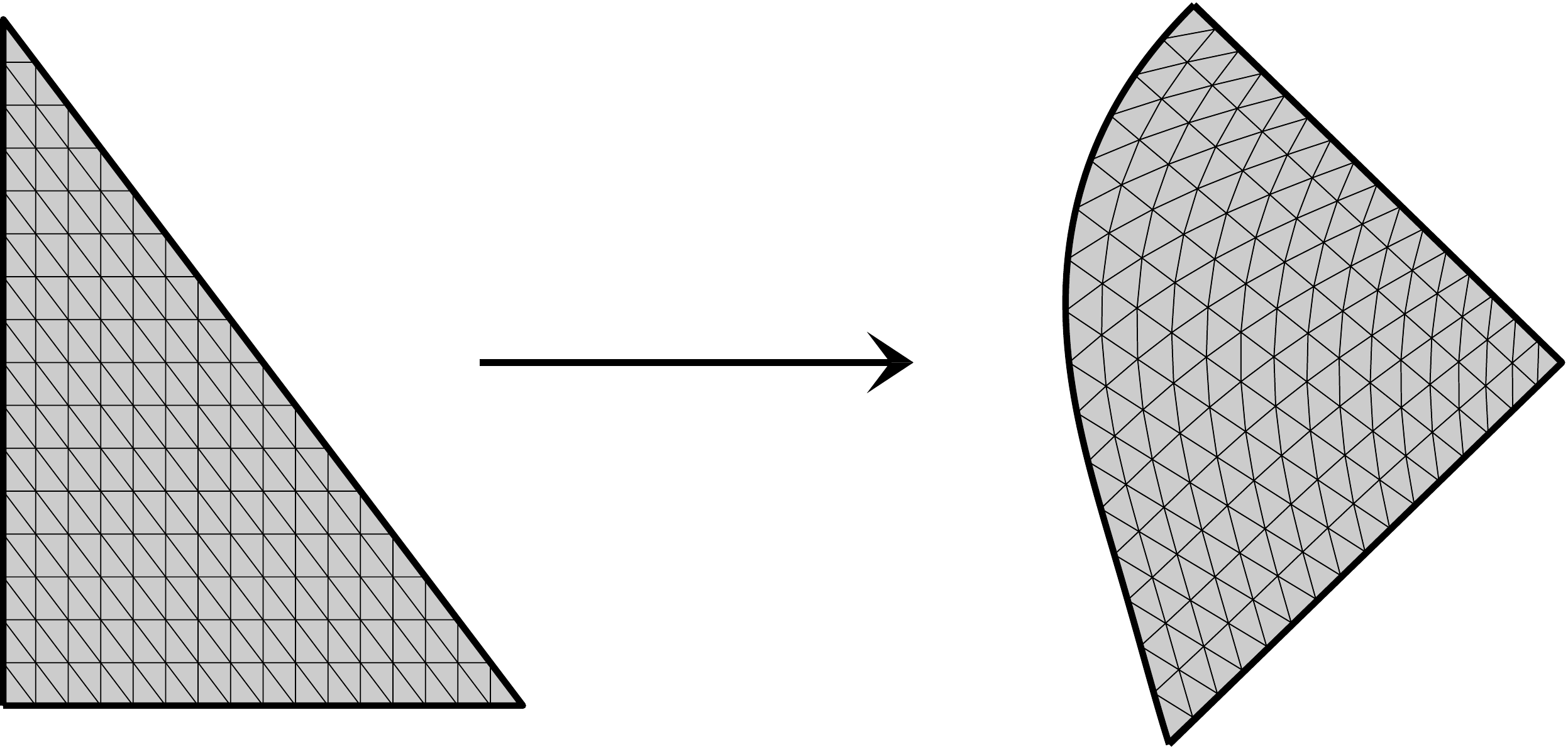}\\
\vspace*{-3.3cm}$\kappa_k$\hspace*{1.5cm}\vspace*{3.3cm}\\[-3ex]
\caption{\label{fig:parameterization}Illustration of the 
diffeomorphim $\kappa_k:\triangle\to\tau_{0,k}$
and the construction of parametric finite elements.}
\end{center}
\end{figure}

The cell functions will be computed by the finite element 
method \cite{B,Brenner}. To construct a triangulation which 
resolves the interface $\partial\omega$ exactly, we use 
\emph{parametric\/} finite elements. To this end, we define a 
macro triangulation with the help of the parametric representation 
of $\partial\omega$ that consists of 28 curved elements based 
on the construction of Zenisek, cf.~\cite{Z}. This macro triangulation
yields a collection of smooth triangular \emph{patches} $\{\tau_{0,k}\}$
and associated diffeomorphisms $\kappa_k:\triangle
\to\tau_{0,k}$ (compare Figure~\ref{fig:parameterization}) 
such that
\begin{equation}\label{eq:parameterization}
  \overline{Y} = \bigcup_{k=1}^{28} \tau_{0,k},
  	\quad \tau_{0,k} = \kappa_k(\triangle), 
	\quad k = 1,2,\ldots,28,
\end{equation}
where $\triangle$ denotes the reference triangle in $\mathbb{R}^2$.
The intersection $\tau_{0,k}\cap \tau_{0,k'}$, $k\not=k'$, of the 
patches $\tau_{0,k}$ and $\tau_{0,k'}$ is either $\emptyset$, or 
a common edge or vertex. Moreover, the diffeomorphisms 
$\kappa_k$ and $\kappa_{k'}$ 
coincide at a common edge except for orientation.

\begin{figure}[hbt]
\begin{center}
\includegraphics[width=6cm,height=6cm]{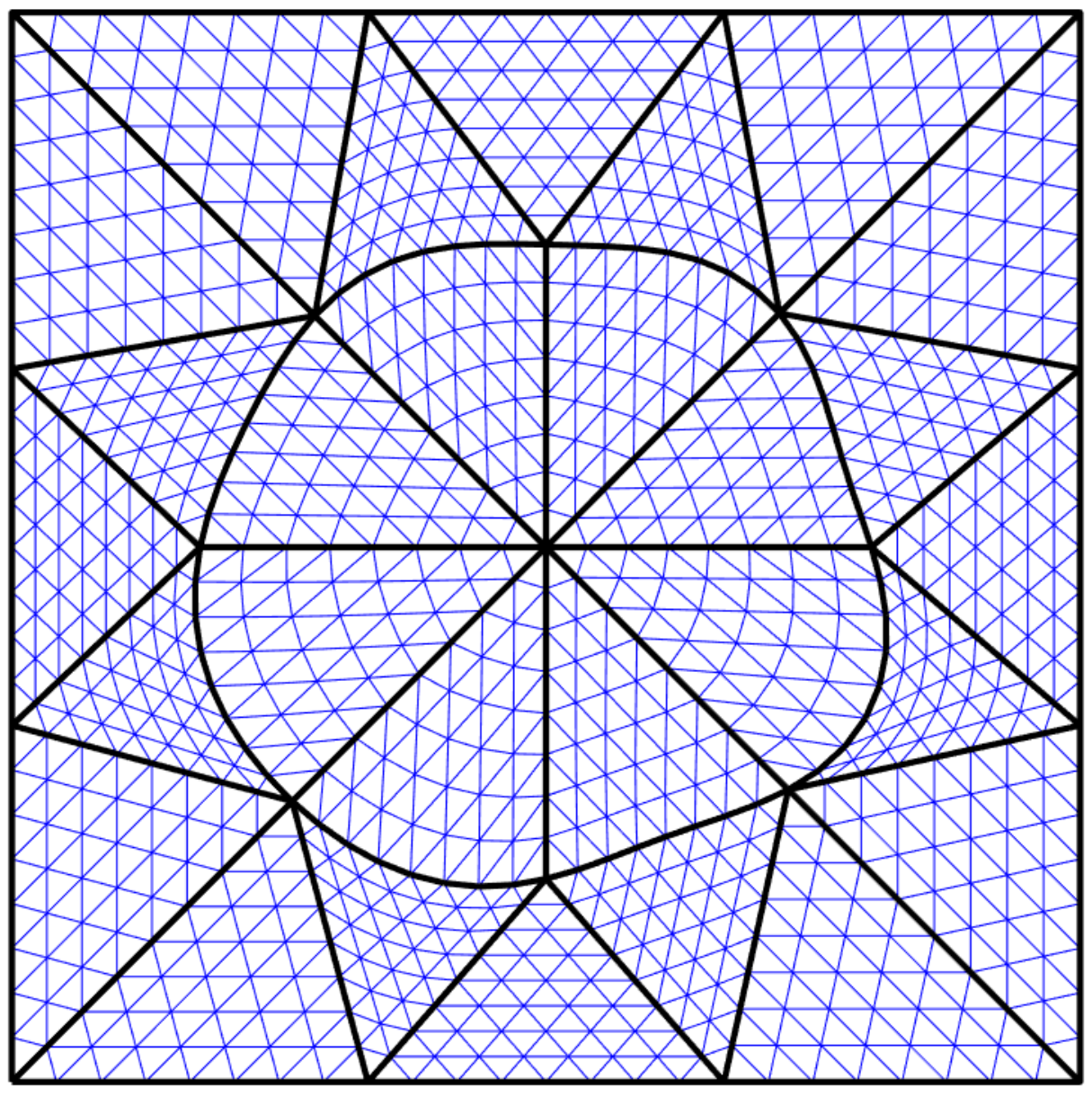}
\caption{\label{fig:mesh}
The macro triangulation consisting of 28 curved elements 
and the resulting mesh of the unit cell ${Y}$, which 
resolves the interface $\partial\omega$.}
\end{center}
\end{figure}

A mesh of level $\ell$ on $Y$ is induced by regular subdivisions
of depth $\ell$ of the reference triangle into four sub-triangles. This 
generates the $4^{\ell}\cdot 28$ triangular elements $\{\tau_{\ell,k}\}$.
An illustration of such a triangulation ($\ell = 3$) is found in 
Figure~\ref{fig:mesh}. On the given triangulation, we employ 
continuous, piecewise linear finite elements to compute the 
cell functions, where the resulting system of linear equations is 
iteratively solved by the conjugate gradient method. Since the meshing 
procedure generates a nested sequence of finite element spaces,
the Bramble-Pasciak-Xu (BPX) preconditioner \cite{BPX} can be 
applied to precondition the iterative solution process. Notice that the 
triangulation resolves the interface and, thus, the convergence 
order of the approximate cell functions is optimal, i.e., second 
order in the mesh size $h$ with respect to the $L^2(Y)$-norm, 
cf.~\cite{LMWZ}.

In the case of a perforated domain, we have to modify our 
finite element implementation correspondingly. In particular, 
the interior of $\omega$ is empty (hence, the macro triangulation
consists of only 20 curved elements) and its boundary $\partial\omega$ 
serves as homogeneous Neumann boundary. These modifications
of the implementation are straightforward, so that we skip 
the details.

For our numerical examples, we choose the expansion 
degree $N=32$ in \eqref{eq:Fourier}, i.e., we consider 65 
design parameters. Moreover, we use roughly $460\,000$ 
finite elements (which corresponds to the refinement level 
$\ell = 7$ of the macro triangulation) for the domain discretization 
in order to solve the state equation for computing the shape functional 
and its gradient. The optimization procedure consists of gradient 
descent method, which is stopped when the $\ell^2$-norm of 
the discrete shape gradient is smaller than $\epsilon = 10^{-5}$.

\subsection{First example}
\begin{figure}[hbt]
\begin{center}
\begin{minipage}{3cm}
\begin{center}
$b_{1,1} = 1.8$\\
\includegraphics[width=3cm,height=3cm]{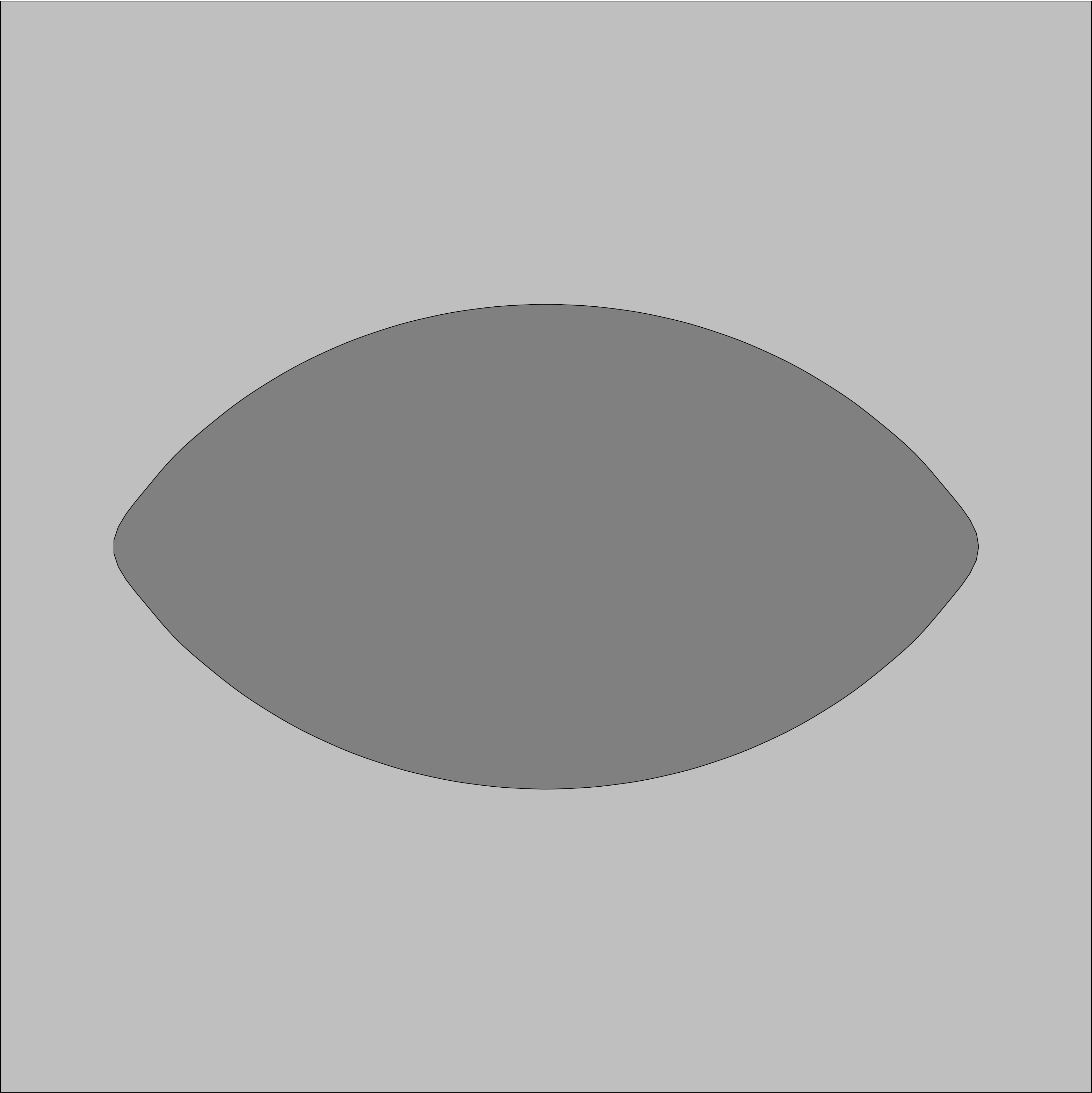}
\end{center}
\end{minipage}
\begin{minipage}{3cm}
\begin{center}
$b_{1,1} = 1.7$\\
\includegraphics[width=3cm,height=3cm]{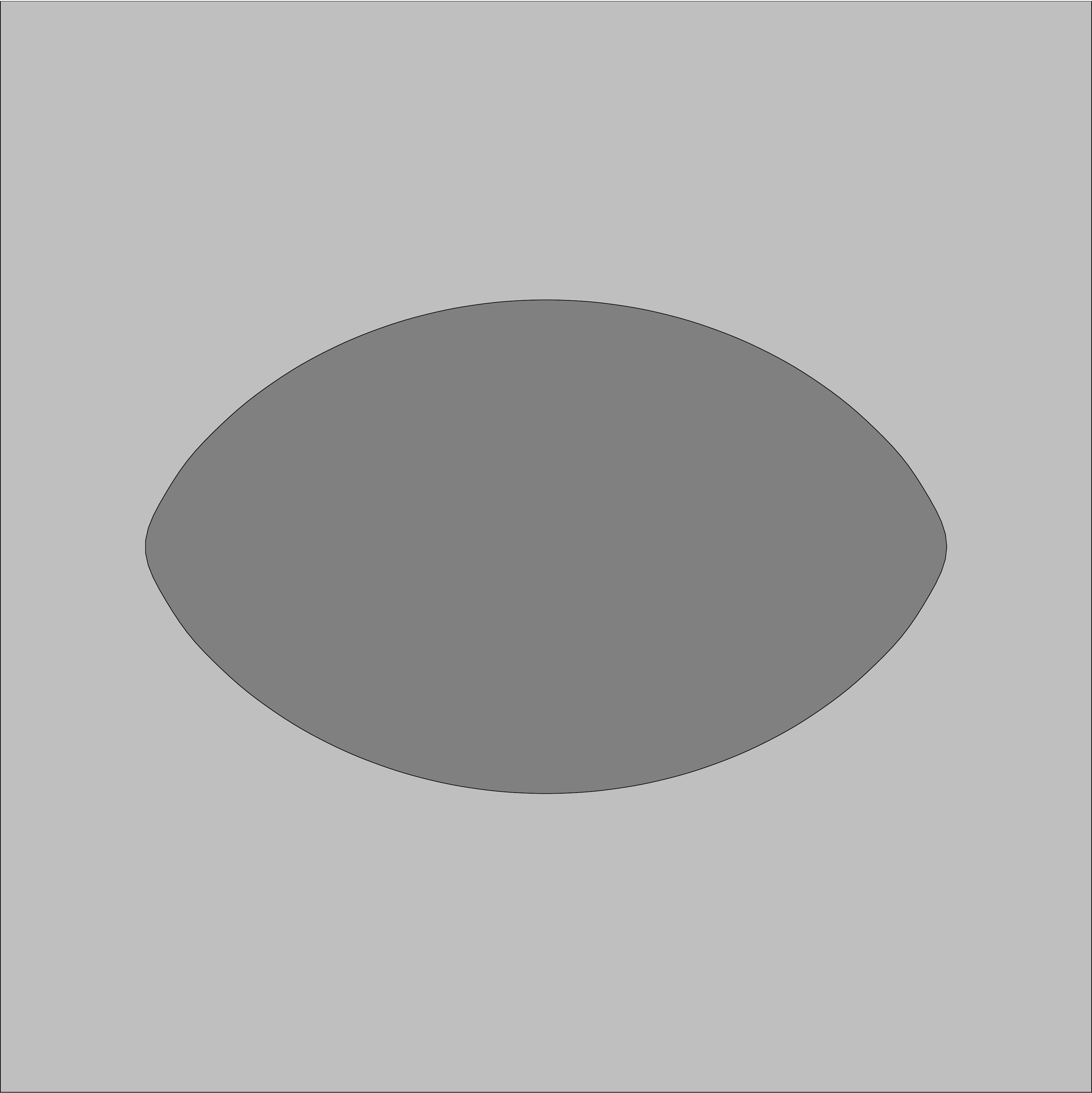}
\end{center}
\end{minipage}
\begin{minipage}{3cm}
\begin{center}
$b_{1,1} = 1.6$\\
\includegraphics[width=3cm,height=3cm]{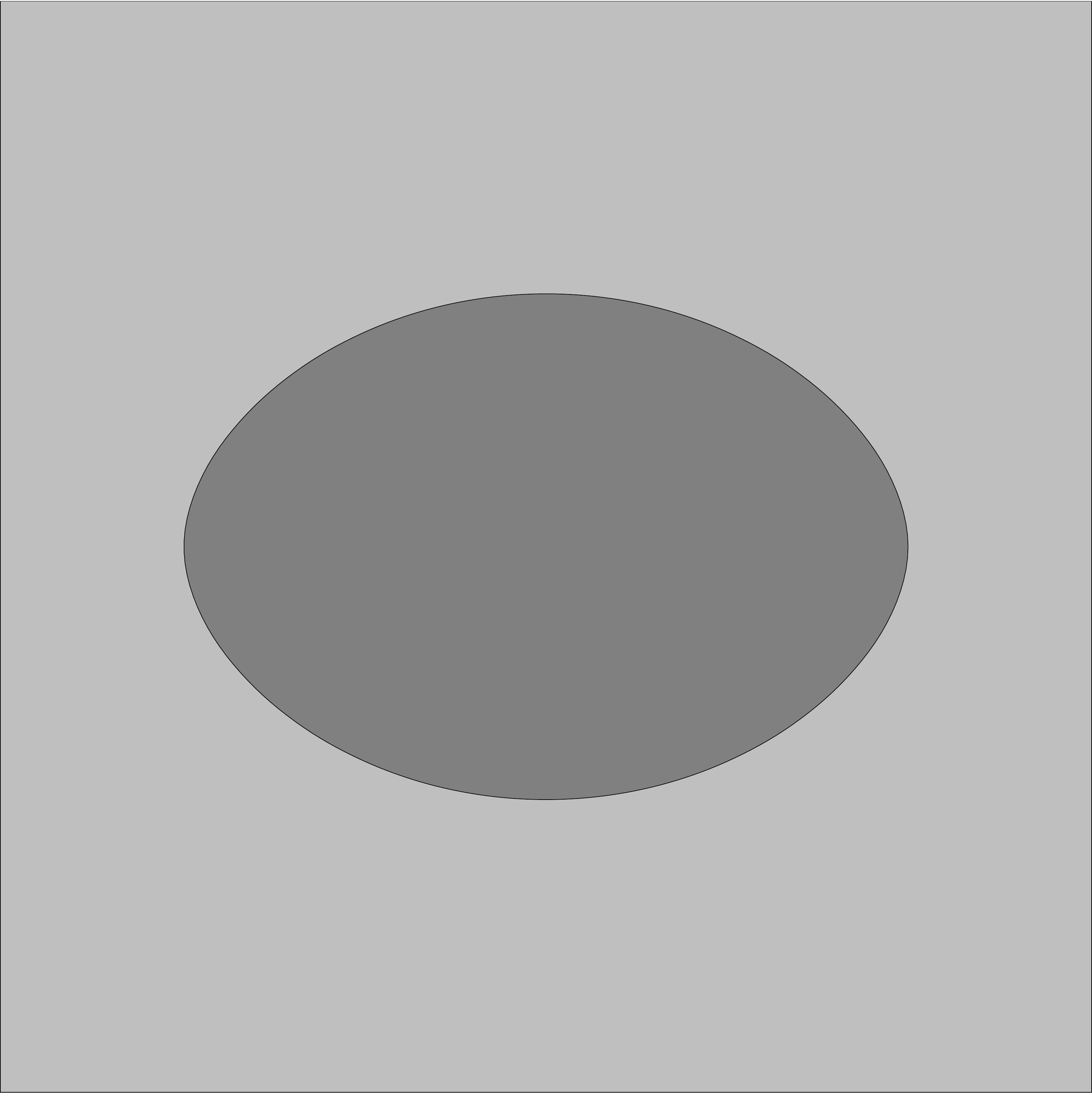}
\end{center}
\end{minipage}
\begin{minipage}{3cm}
\begin{center}
$b_{1,1} = 1.5$\\
\includegraphics[width=3cm,height=3cm]{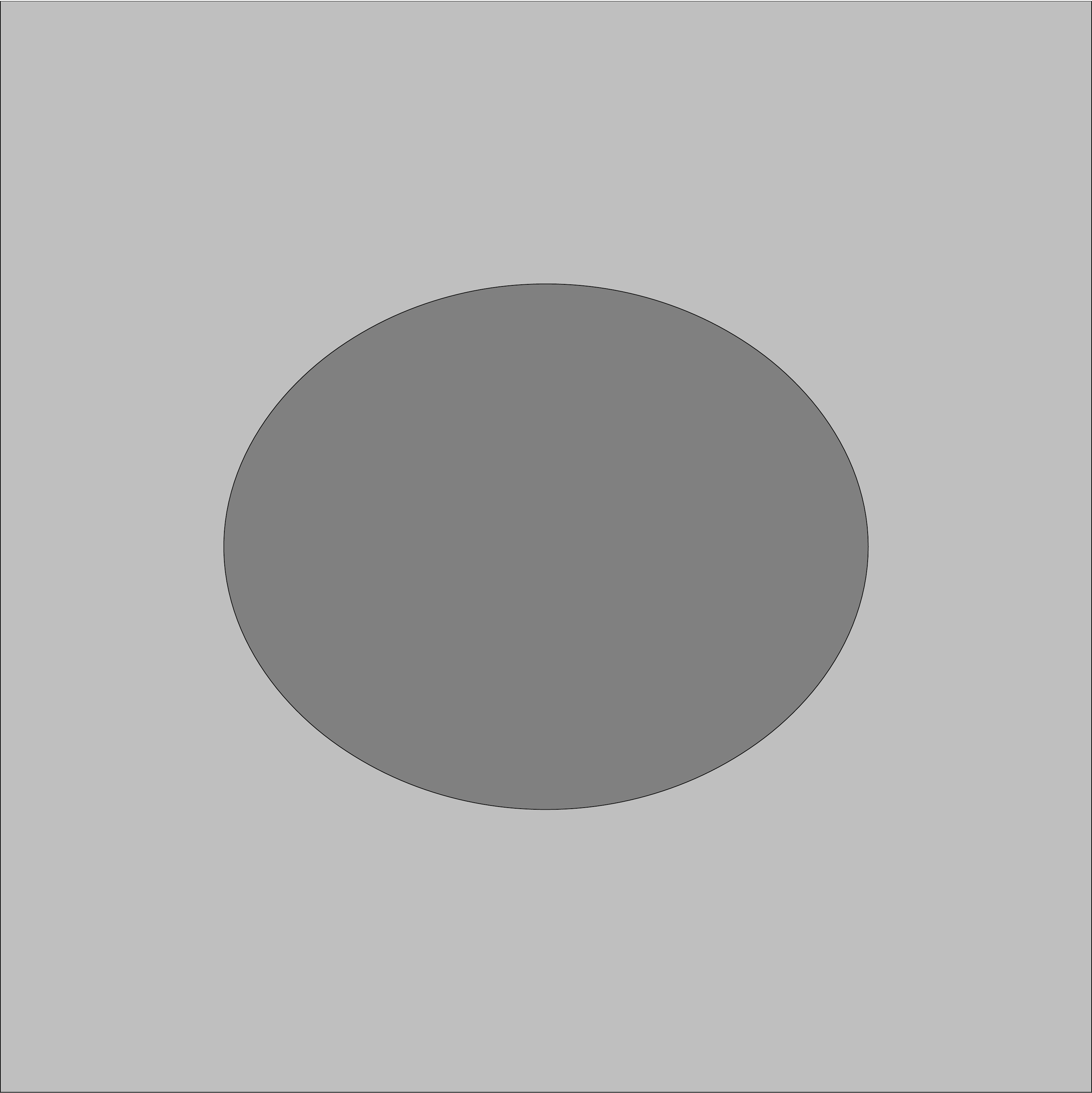}
\end{center}
\end{minipage}\\[1ex]
\begin{minipage}{3cm}
\begin{center}
$b_{1,1} = 1.4$\\
\includegraphics[width=3cm,height=3cm]{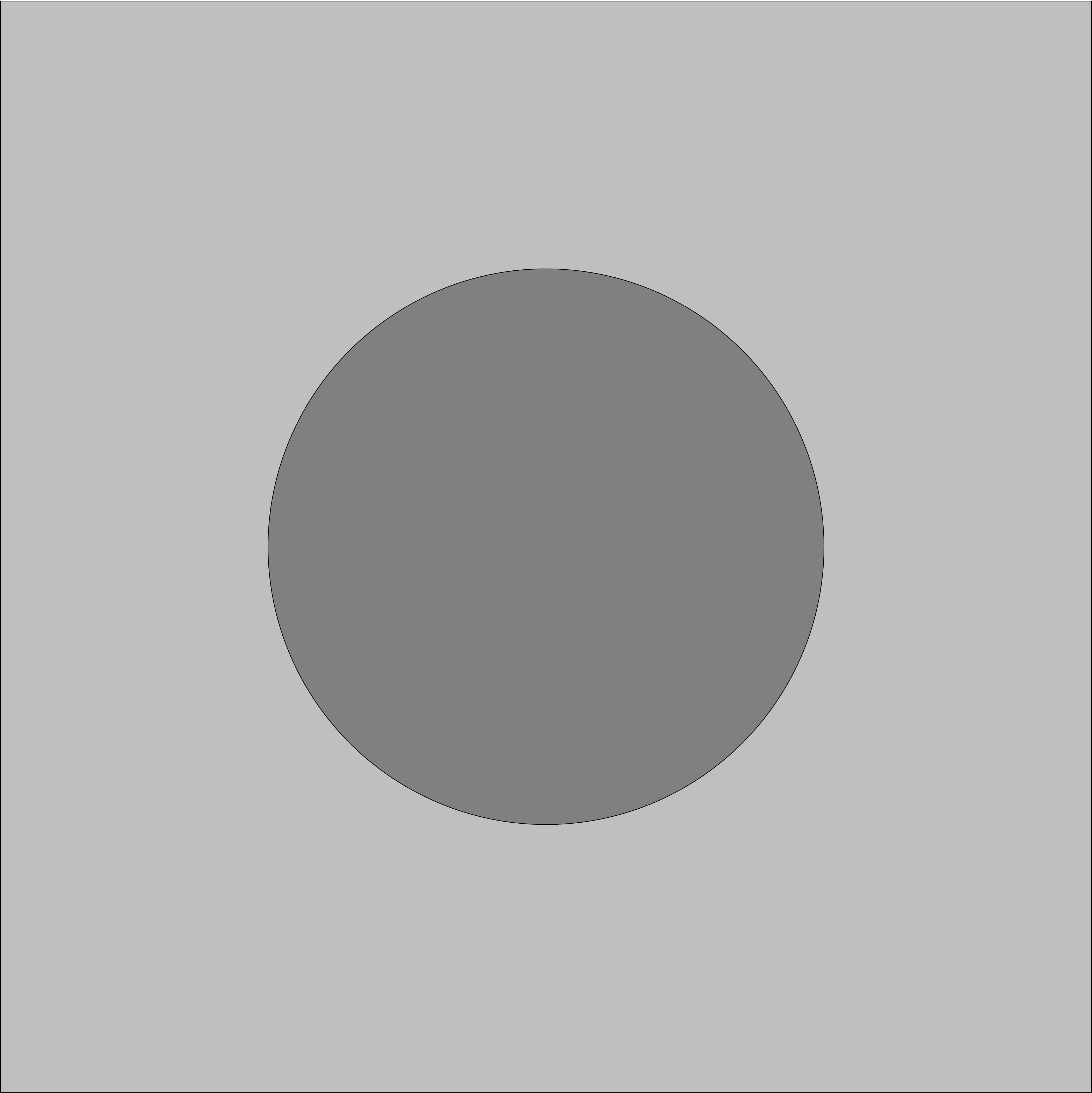}
\end{center}
\end{minipage}
\begin{minipage}{3cm}
\begin{center}
$b_{1,1} = 1.3$\\
\includegraphics[width=3cm,height=3cm]{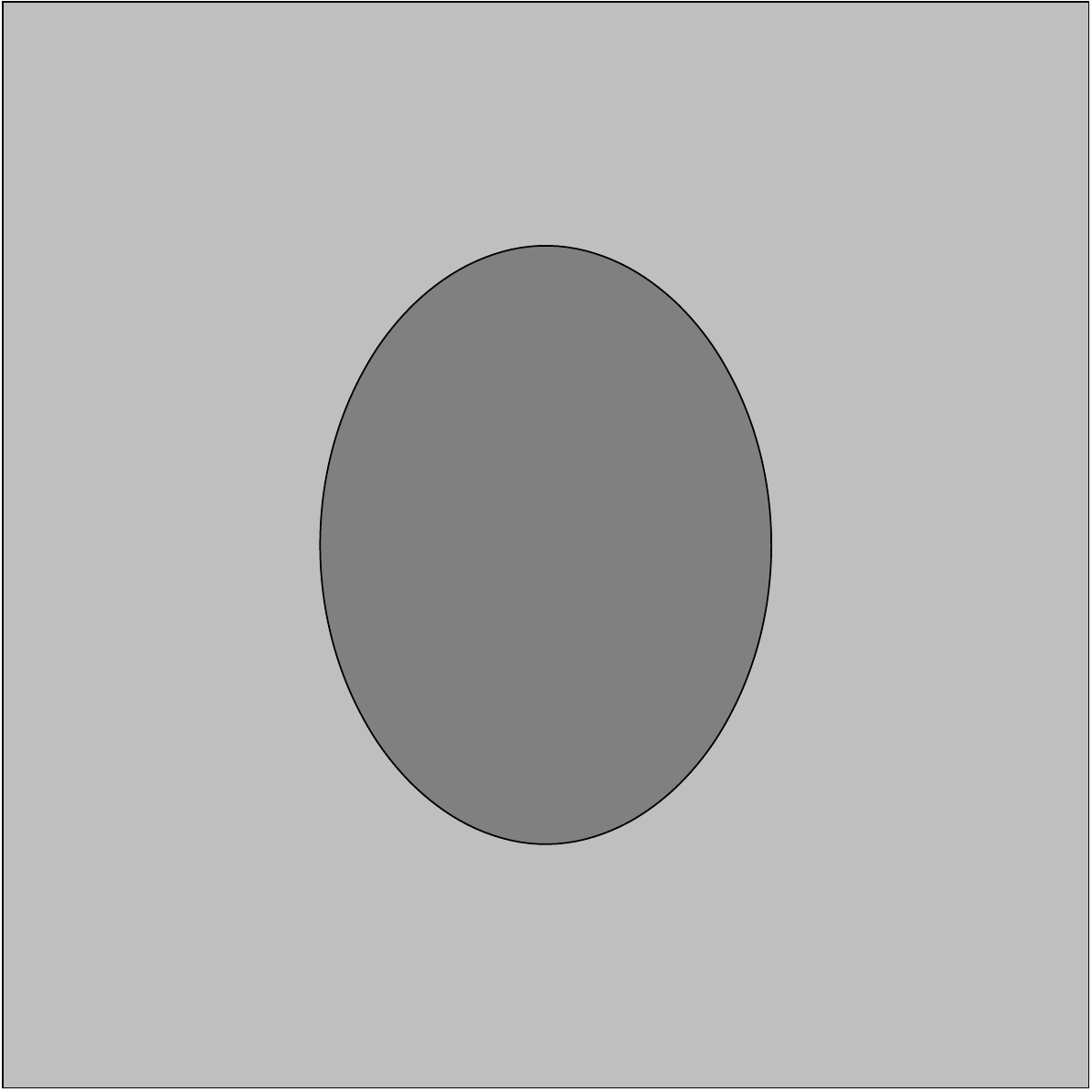}
\end{center}
\end{minipage}
\begin{minipage}{3cm}
\begin{center}
$b_{1,1} = 1.2$\\
\includegraphics[width=3cm,height=3cm]{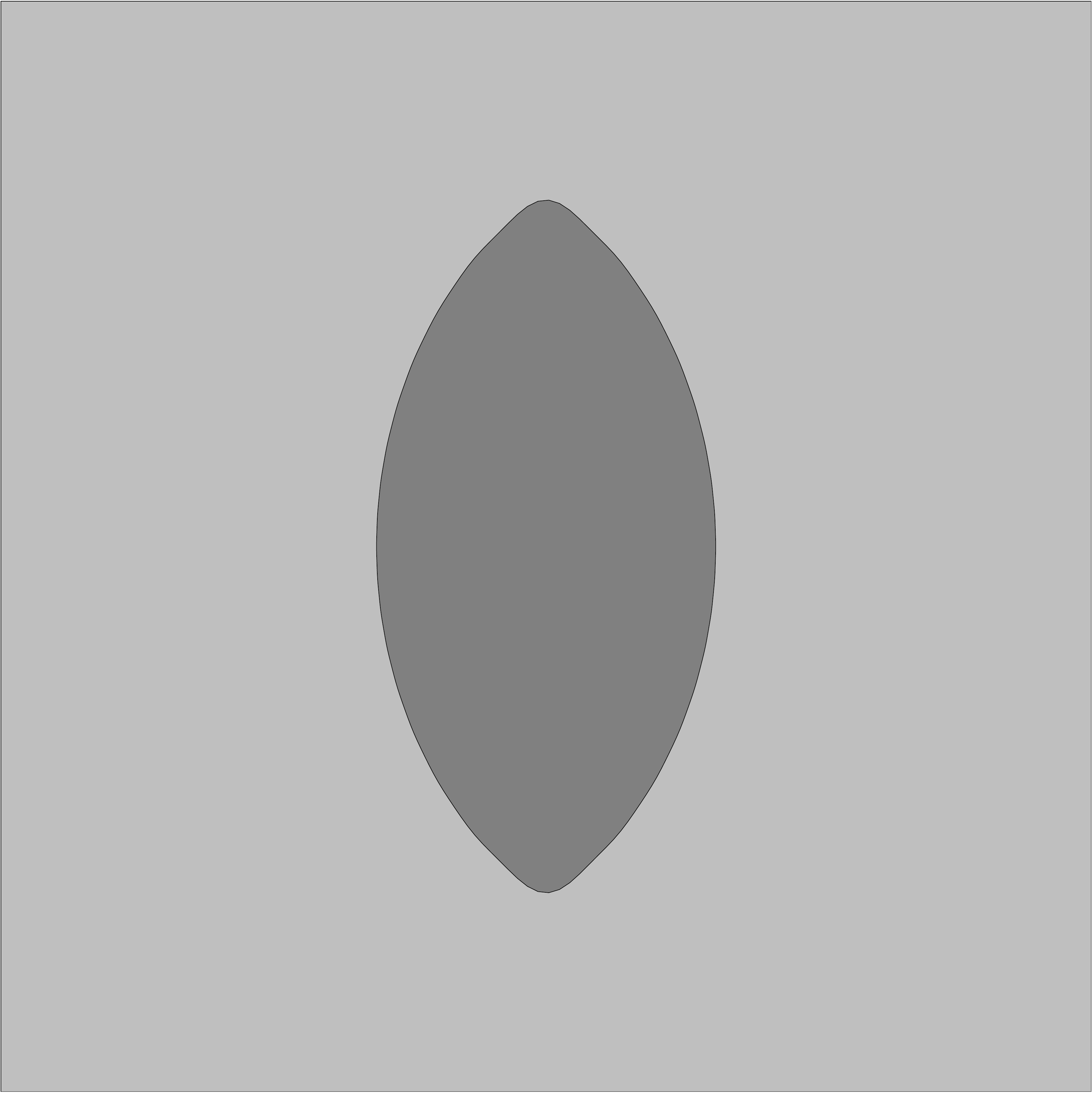}
\end{center}
\end{minipage}
\caption{\label{fig:shapes1}Optimal shapes for the desired 
effective tensor $\boldsymbol{B}_1$ in case of different values 
of $b_{1,1}$ when the circle of radius $1/4$ is used as initial guess.}
\end{center}
\end{figure}

For our first computations, we choose constant coefficient
functions $\sigma_1\equiv 1$ and $\sigma_2\equiv 10$. 
The desired effective tensor in \eqref{eq:functional} is 
\[
\boldsymbol{B}_1 = \begin{bmatrix} b_{1,1} & 0 \\ 0 & 1.4 \end{bmatrix},
\]
where $b_{1,1}$ varies from $1.1$ to $1.8$ with step size
$0.1$. Starting with the circle of radius $1/4$ as initial guess, 
we obtain the optimal shapes found in Figure~\ref{fig:shapes1}.
Note that the desired effective tensor is achieved, i.e., the 
shape functional \eqref{eq:functional} is zero in the computed 
shapes. We especially see that the alignment of the computed 
shape reflects the anisotropy of the effective tensor.

\subsection{Second example}
\begin{figure}[hbt]
\begin{center}
\begin{minipage}{3cm}
\begin{center}
$b_{1,2} = -0.10$\\
\includegraphics[width=3cm,height=3cm]{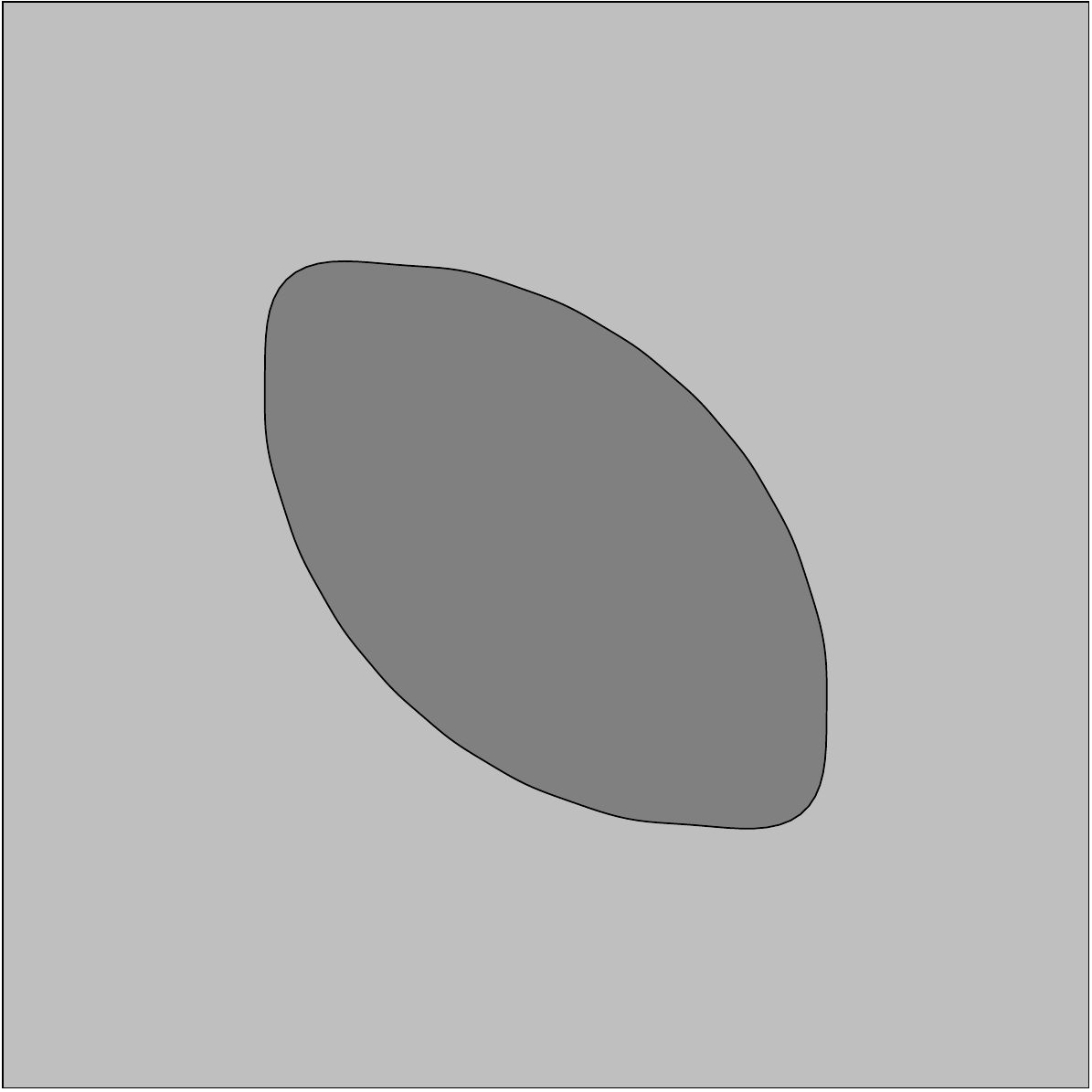}
\end{center}
\end{minipage}
\begin{minipage}{3cm}
\begin{center}
$b_{1,2} = -0.05$\\
\includegraphics[width=3cm,height=3cm]{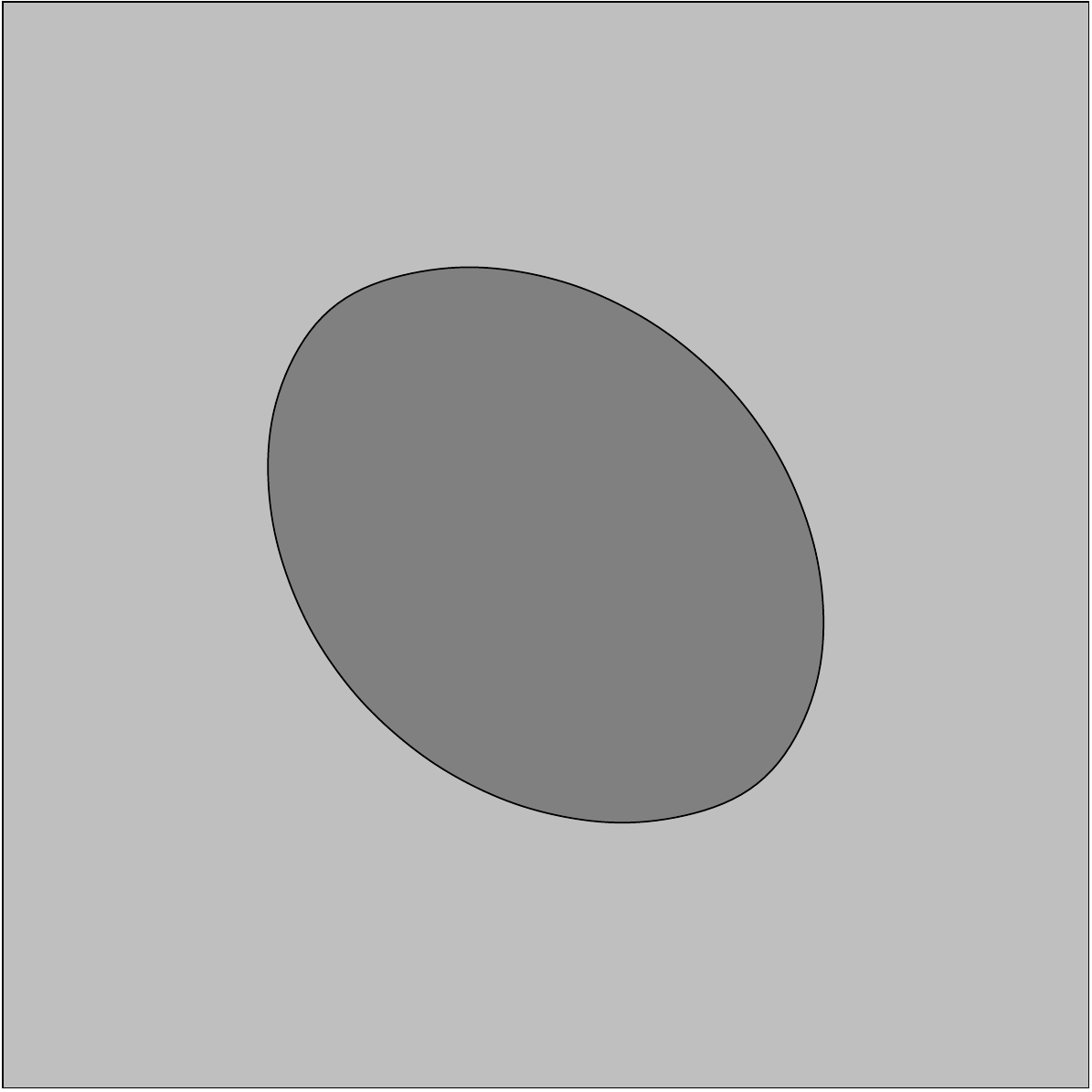}
\end{center}
\end{minipage}
\begin{minipage}{3cm}
\begin{center}
$b_{1,2} = 0.05$\\
\includegraphics[width=3cm,height=3cm]{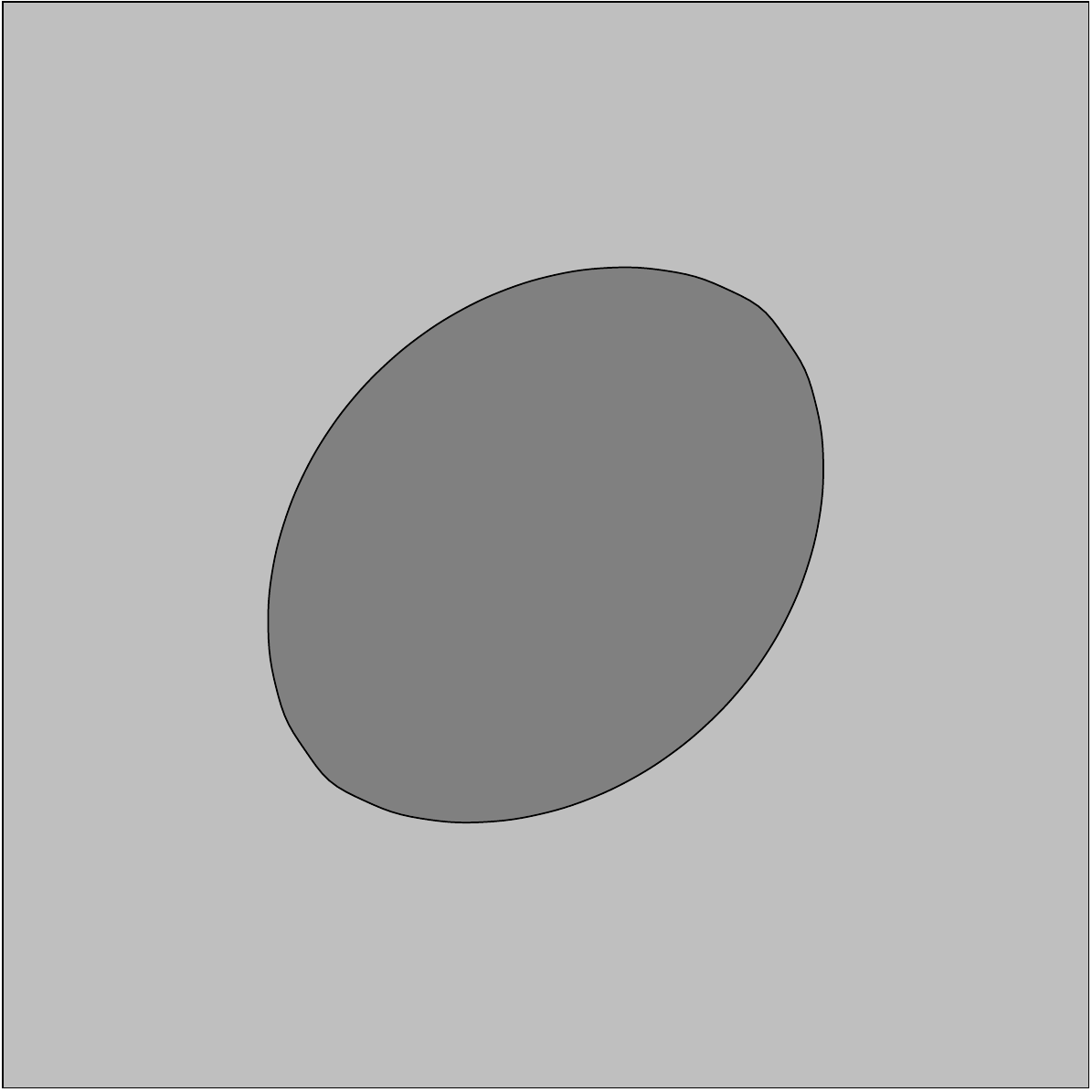}
\end{center}
\end{minipage}
\begin{minipage}{3cm}
\begin{center}
$b_{1,2} = 0.10$\\
\includegraphics[width=3cm,height=3cm]{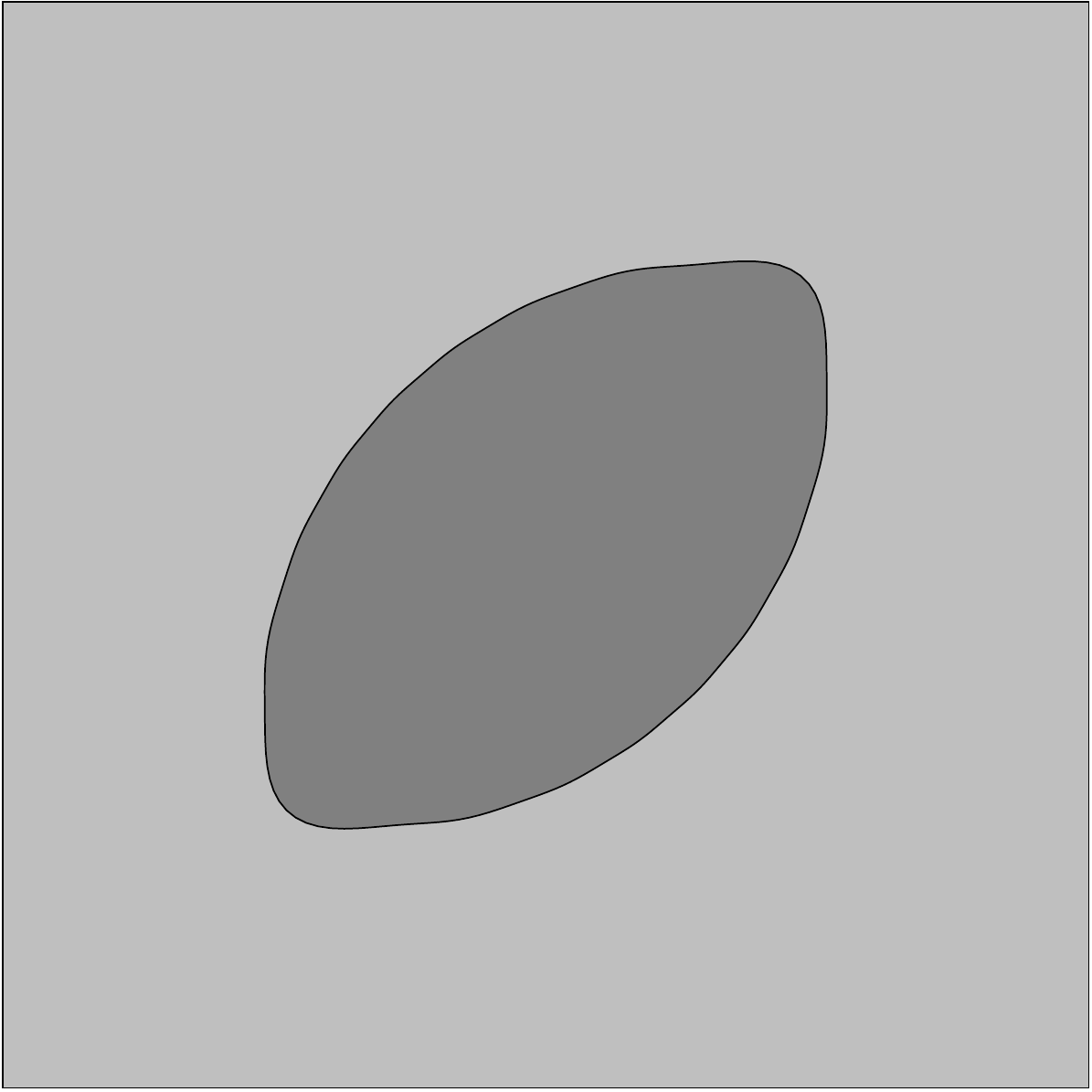}
\end{center}
\end{minipage}
\caption{\label{fig:shapes2}Optimal shapes for the desired 
effective tensor $\boldsymbol{B}_2$ in case of different values 
of $b_{1,2} = b_{2,1}$ when the circle of radius $1/4$ is used 
as initial guess.}
\end{center}
\end{figure}

We shall study the effect of the off-diagonal terms
in the desired effective tensor. We thus consider
\[
\boldsymbol{B}_2 = \begin{bmatrix} 1.4 & b_{1,2} \\ b_{2,1} & 1.4 \end{bmatrix}
\]
with $b_{1,2} = b_{2,1}$ chosen to be equal to $\pm 0.05$ 
and $\pm 0.1$. The results are found in Figure~\ref{fig:shapes2}
in the order $-0.10$, $-0.05$, $0.05$, $0.10$ (from left to right) for the
values of $b_{1,2} = b_{2,1}$. It is seen that the shape is oriented 
north-west in case of a negative sign and north-east in case of 
a positive sign. Notice that we obtain the circle found in Figure
\ref{fig:shapes1} in the situation $b_{1,2} = b_{2,1} = 0$.

\subsection{Third example}
\begin{figure}[hbt]
\begin{center}
\includegraphics[width=3cm,height=3cm]{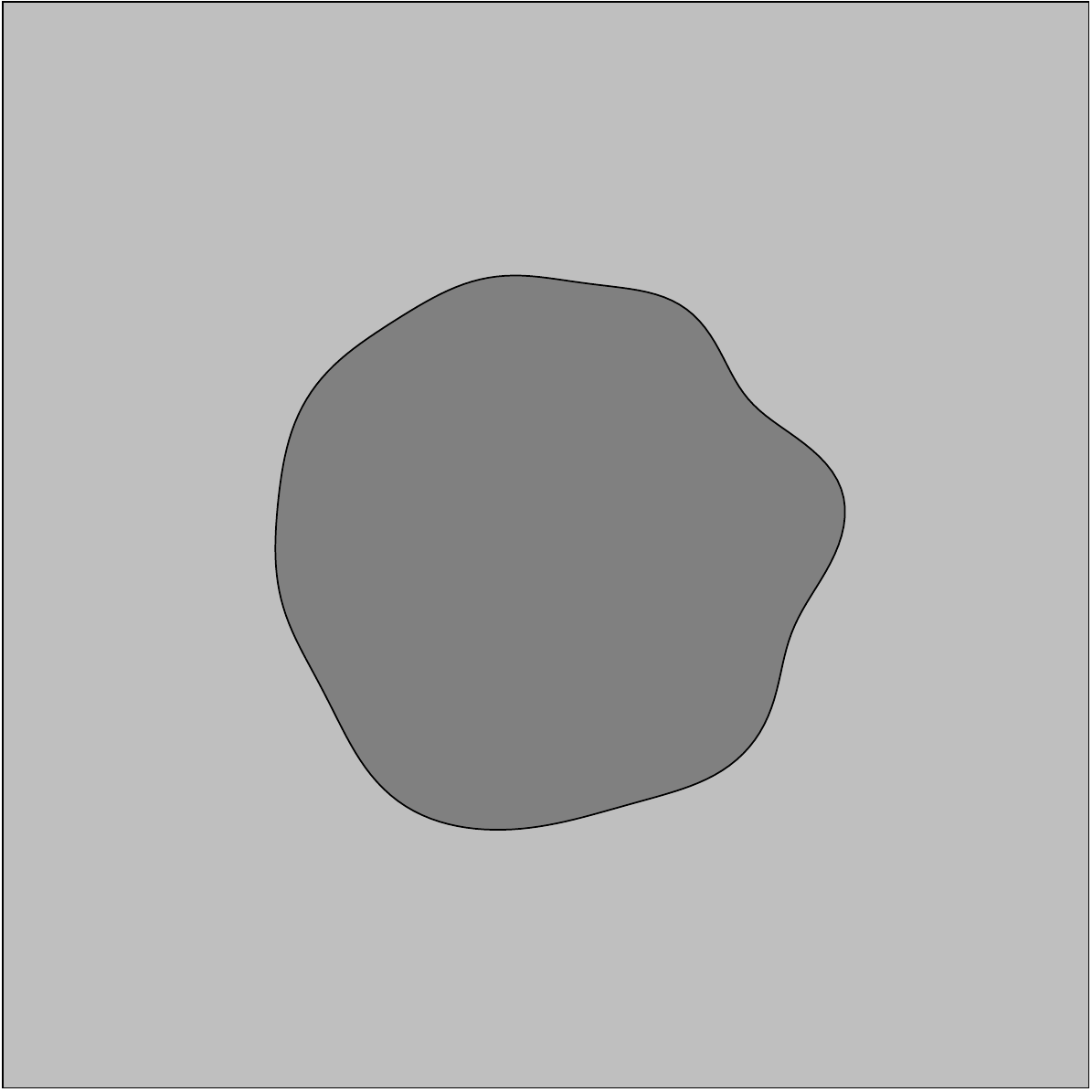}
\includegraphics[width=3cm,height=3cm]{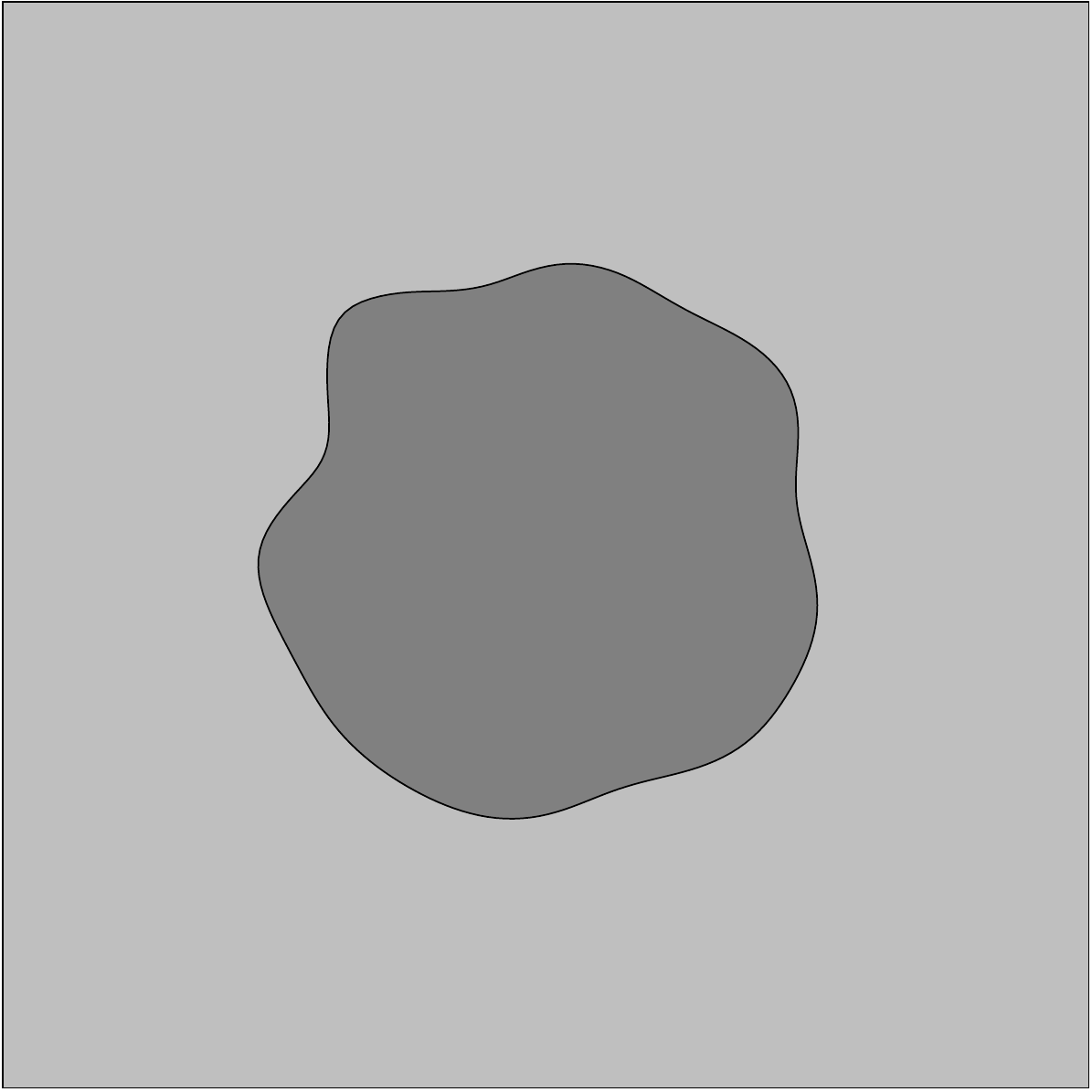}
\includegraphics[width=3cm,height=3cm]{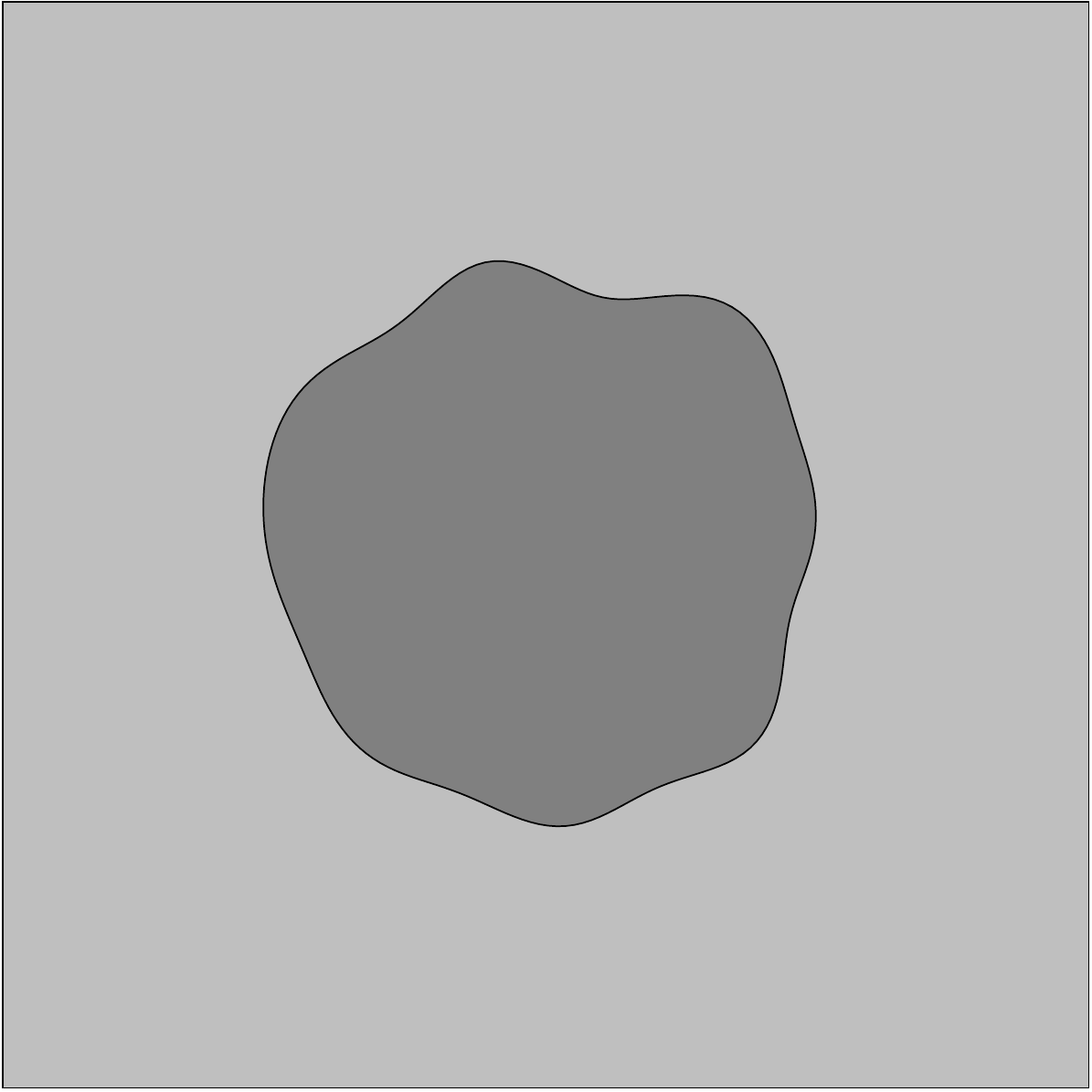}
\includegraphics[width=3cm,height=3cm]{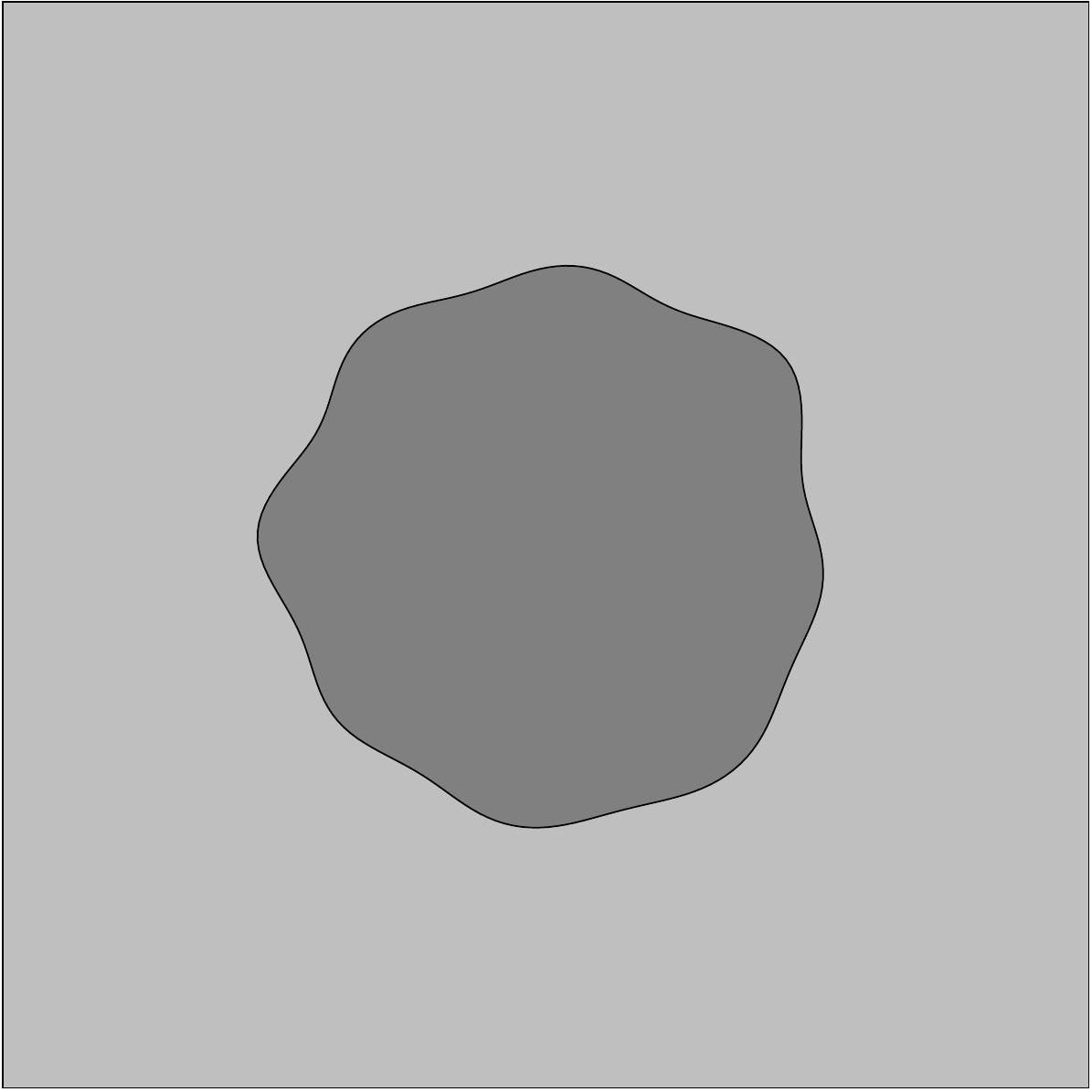}
\caption{\label{fig:shapes3}Optimal shapes for the desired 
effective tensor $\boldsymbol{B}_3$ when starting with a 
randomly perturbed circle as initial guess.}
\end{center}
\end{figure}

In our next test, we shall show that the solution for $\omega$ 
is non-unique. To this end, we choose a randomly perturbed 
circle of radius $1/4$ as initial guess and try to construct a 
microstructure that has the (isotropic) effective tensor
\begin{equation}\label{eq:B3}
\boldsymbol{B}_3 = \begin{bmatrix} 1.4 & 0 \\ 0 & 1.4 \end{bmatrix}.
\end{equation}
In Figure~\ref{fig:shapes3}, we see the different shapes we
get from the minimization of the shape functional \eqref{eq:functional},
all of them resulting in the desired effective tensor $\boldsymbol{B}_3$.
Notice that we obtain a circle if we would start with a circle 
as it can be seen in the fifth plot in Figure~\ref{fig:shapes1}.

\subsection{Fourth example}
\begin{figure}[hbt]
\begin{center}
\includegraphics[width=3cm,height=3cm]{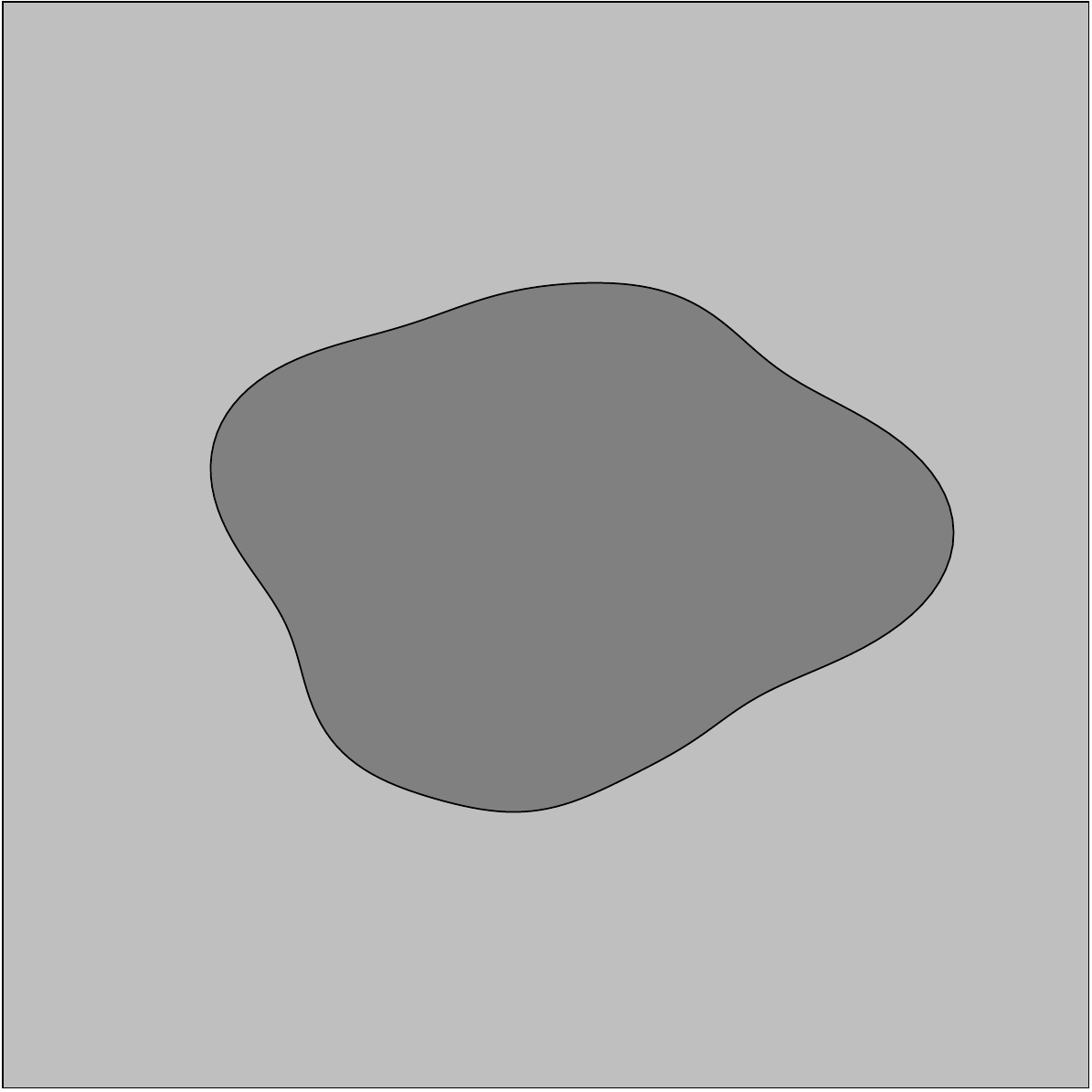}
\includegraphics[width=3cm,height=3cm]{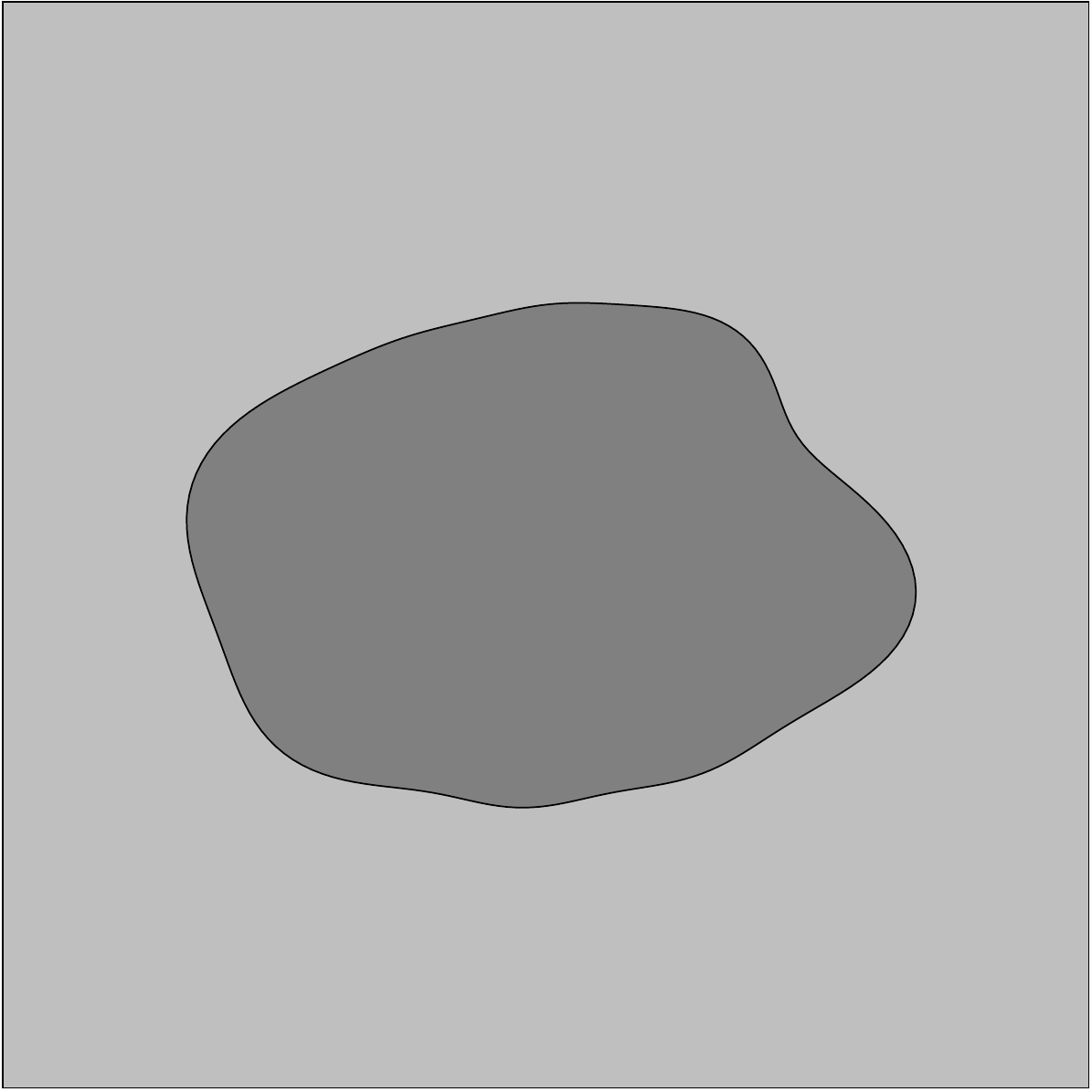}
\includegraphics[width=3cm,height=3cm]{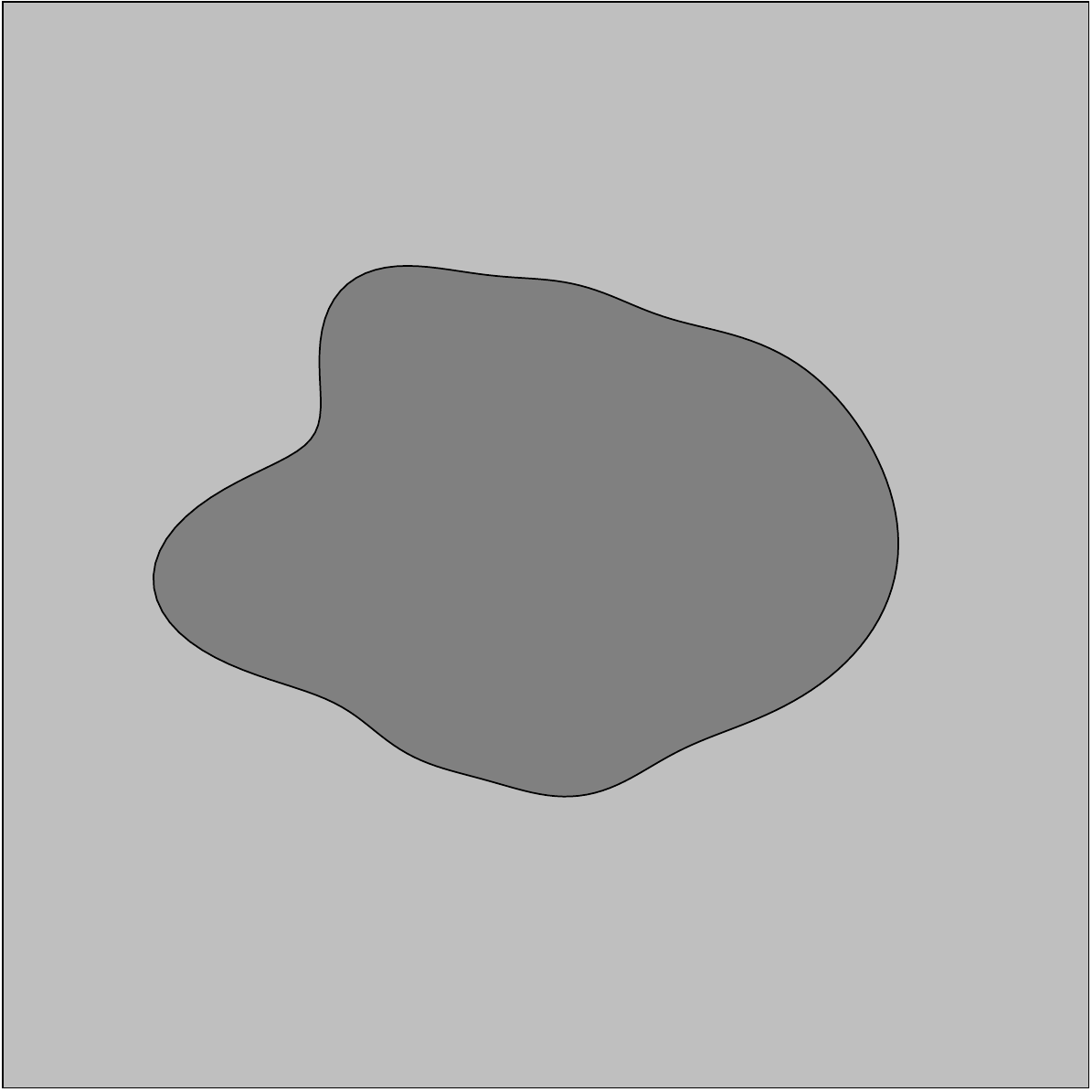}
\includegraphics[width=3cm,height=3cm]{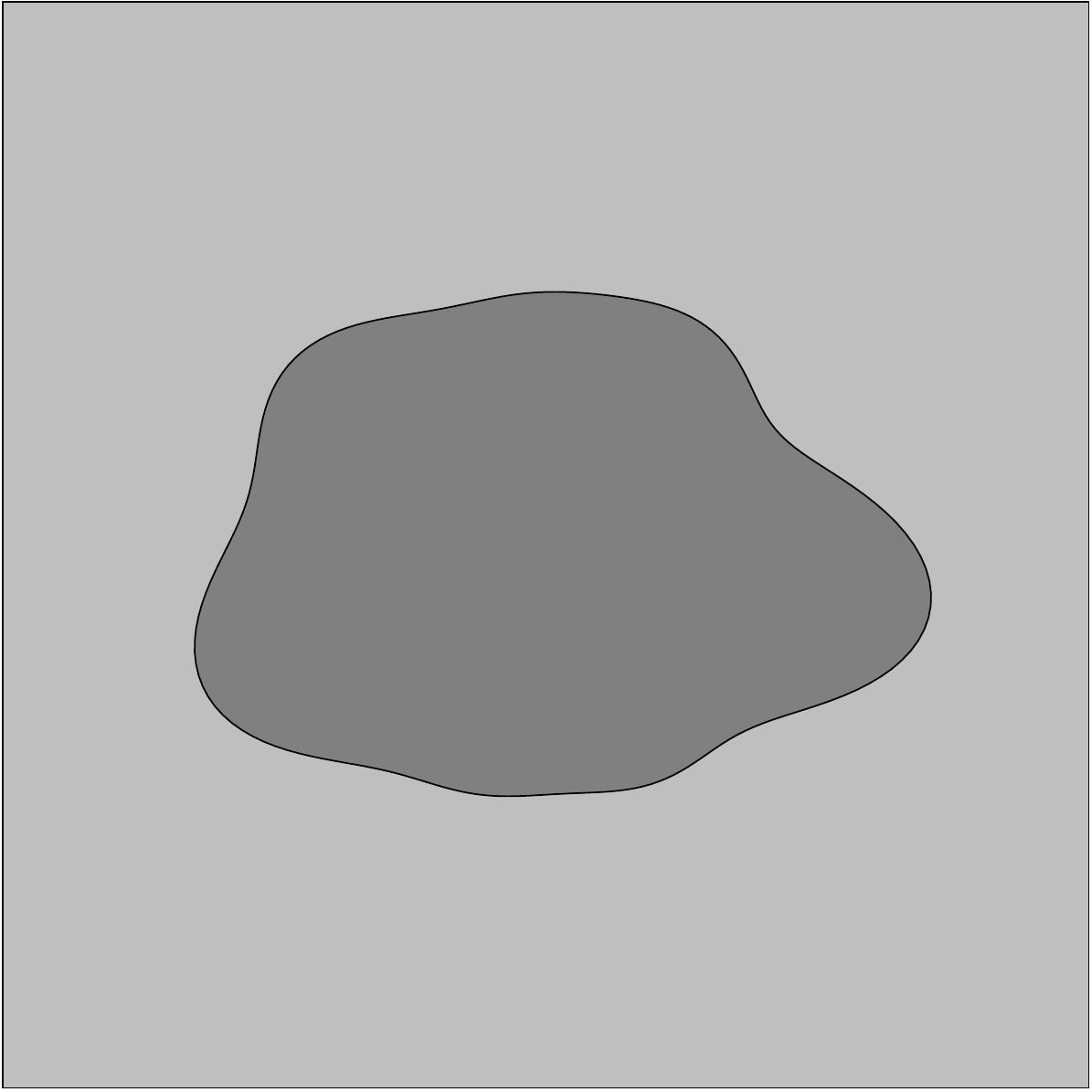}
\caption{\label{fig:shapes4}Optimal shapes for the desired 
effective tensor $\boldsymbol{B}_4$ when starting with a 
randomly perturbed circle as initial guess.}
\end{center}
\end{figure}

Likewise, we observe the same non-uniqueness if the 
desired effective tensor is anisotropic. Choosing
\[
\boldsymbol{B}_4 = \begin{bmatrix} 1.6 & 0 \\ 0 & 1.4 \end{bmatrix}
\]
and starting again by a randomly perturbed circle 
of radius $1/4$, we obtain the shapes found in 
Figure~\ref{fig:shapes4}. All these shapes yield 
again the effective tensor $\boldsymbol{B}_4$.
In case of starting with the circle of radius $1/4$, 
we get the shape seen in the third plot in 
Figure~\ref{fig:shapes1}.

\subsection{Fifth example}
\begin{figure}[hbt]
\begin{center}
\includegraphics[width=3cm,height=3cm]{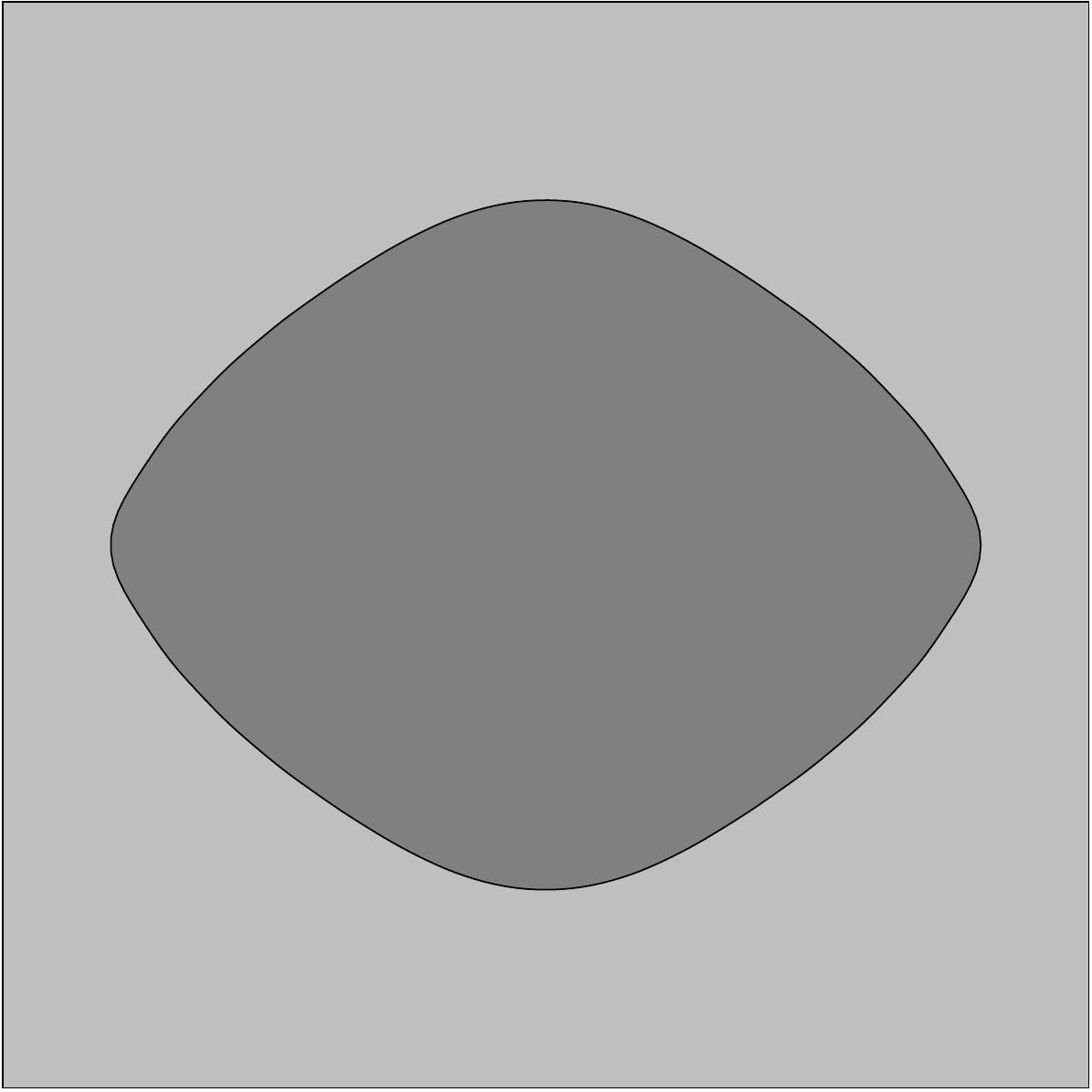}
\includegraphics[width=3cm,height=3cm]{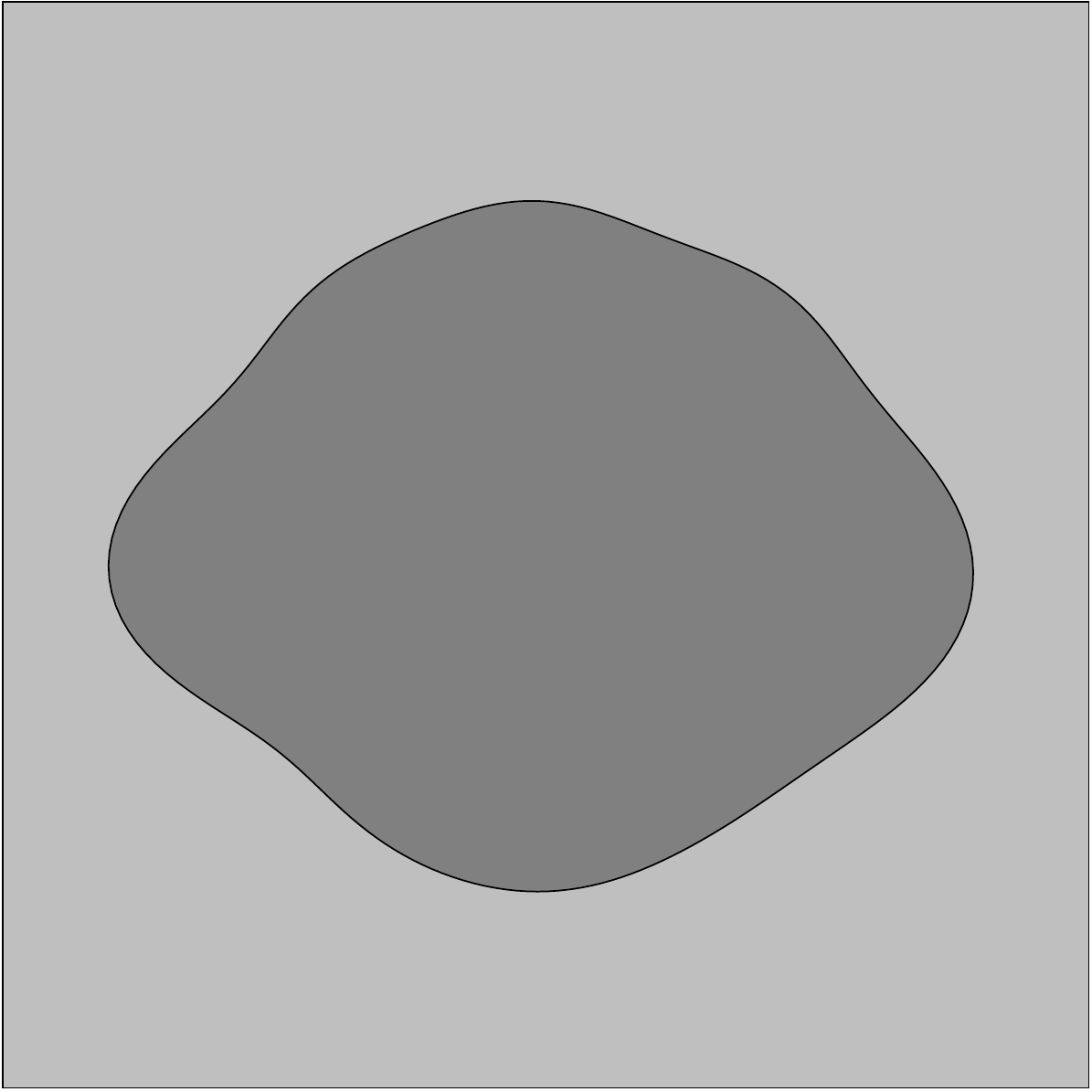}
\includegraphics[width=3cm,height=3cm]{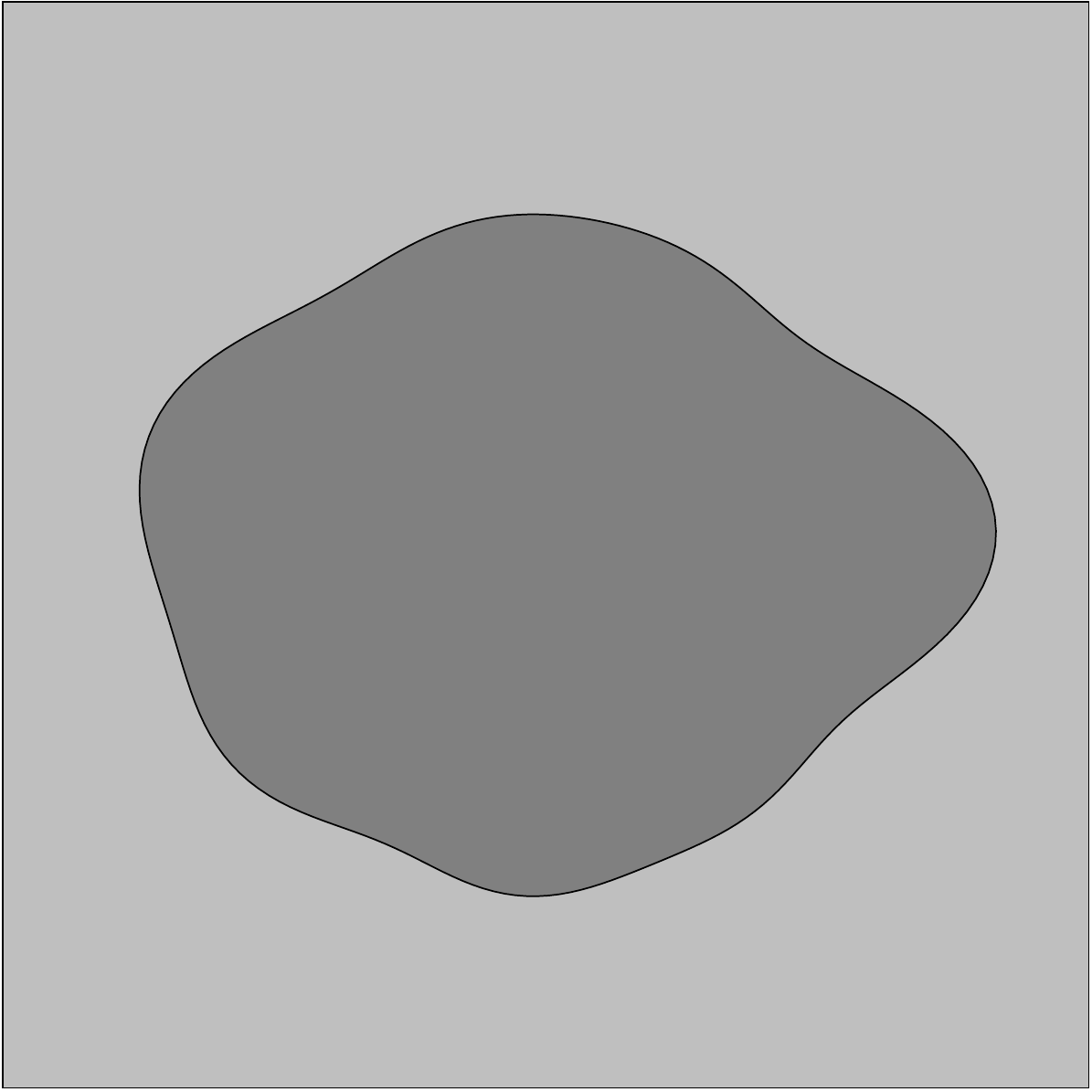}
\includegraphics[width=3cm,height=3cm]{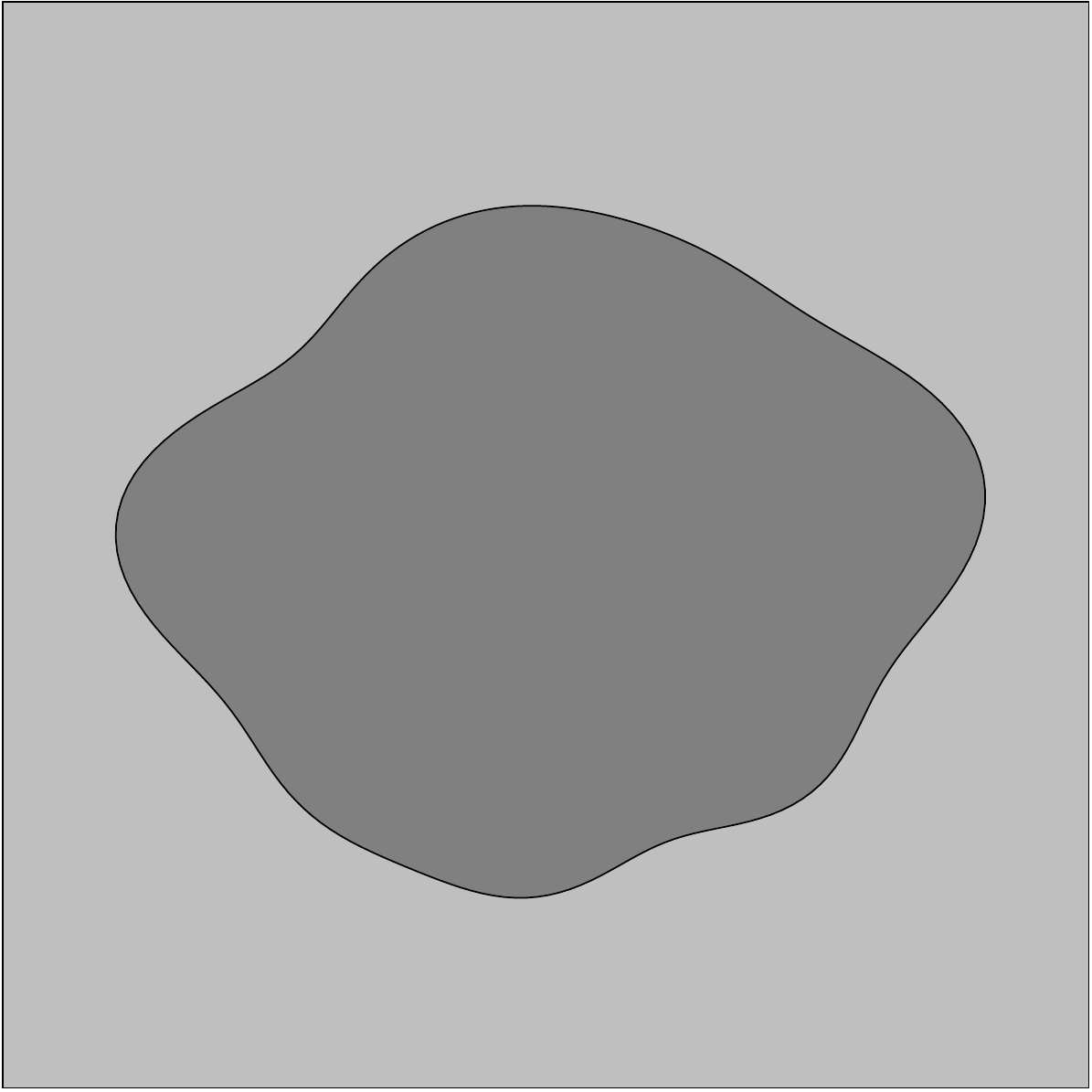}
\caption{\label{fig:shapes5}Different shapes which all generate 
the desired effective tensor $\boldsymbol{B}_3$ in case 
of a non-constant coefficient function $\sigma_2$.}
\end{center}
\end{figure}

We shall also consider the situation that $\sigma_1$ 
and $\sigma_2$ are smooth functions. To this end, we 
consider $\sigma_1\equiv 1$ to be constant but 
\[
  \sigma_2(x,y) = 5\Bigg(\frac{11}{10}+\cos(2\pi x)+4\bigg(y-\frac{1}{2}\bigg)^2\Bigg).
\]
For the desired (isotropic) effective tensor
$\boldsymbol{B}_3$ from \eqref{eq:B3} and the circle with 
radius $1/4$ as initial guess, we obtain the optimal shape
found in the outermost left plot of Figure~\ref{fig:shapes5}.
If we randomly perturb the initial circle, then we obtain 
optimal shapes which are different, compare the other 
plots of Figure~\ref{fig:shapes5} for some results.

\subsection{Sixth example}
\begin{figure}[hbt]
\begin{center}
\includegraphics[width=3cm,height=3cm]{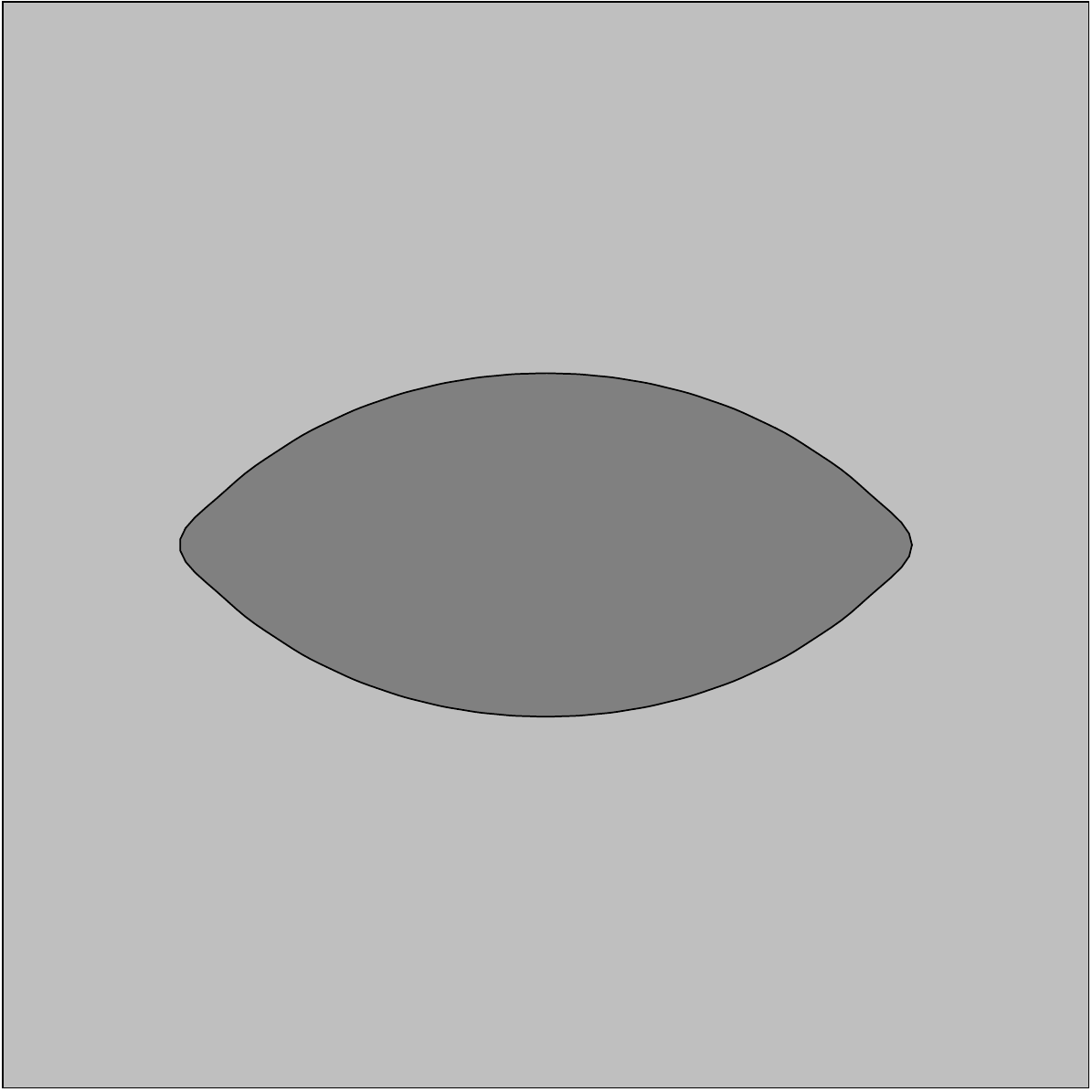}
\includegraphics[width=3cm,height=3cm]{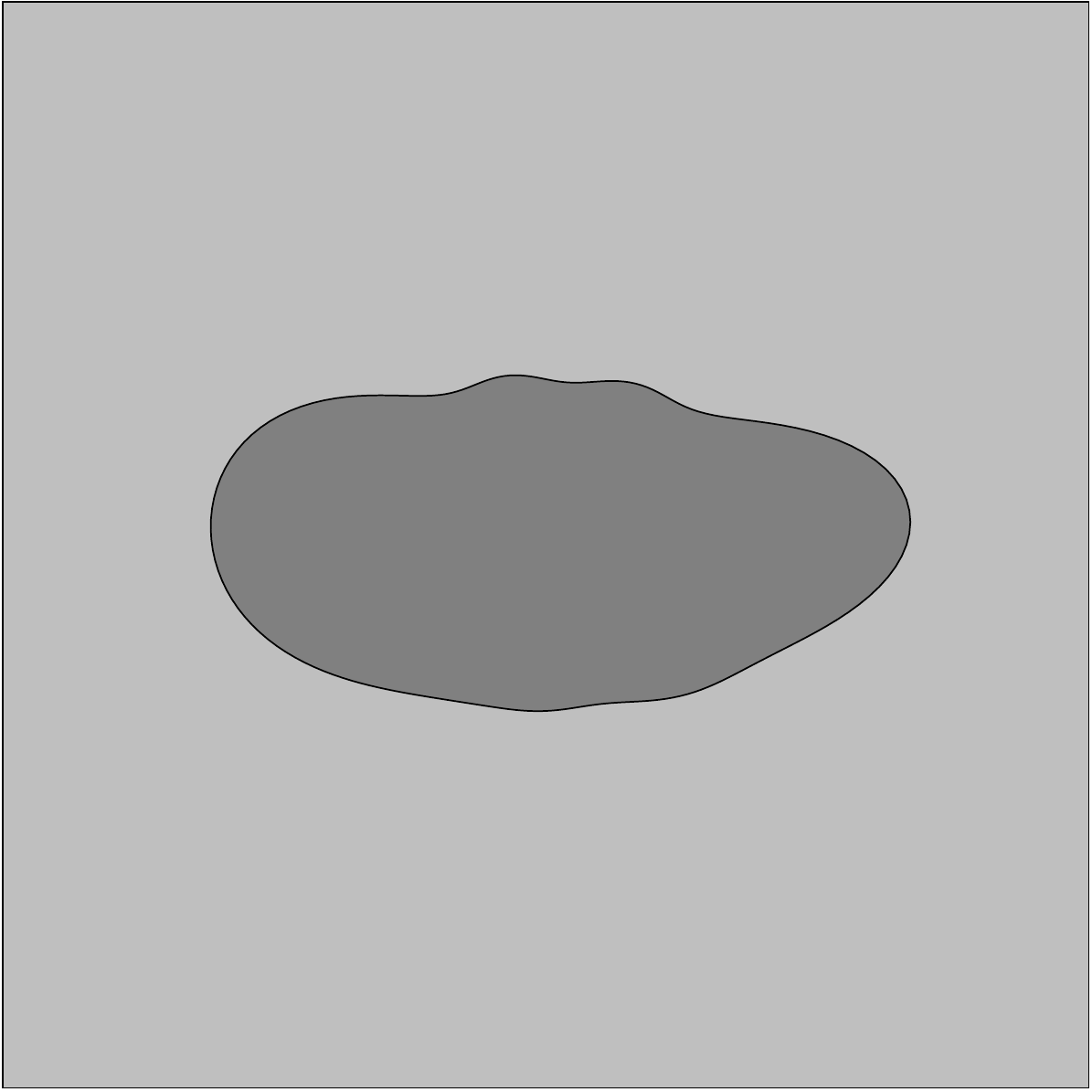}
\includegraphics[width=3cm,height=3cm]{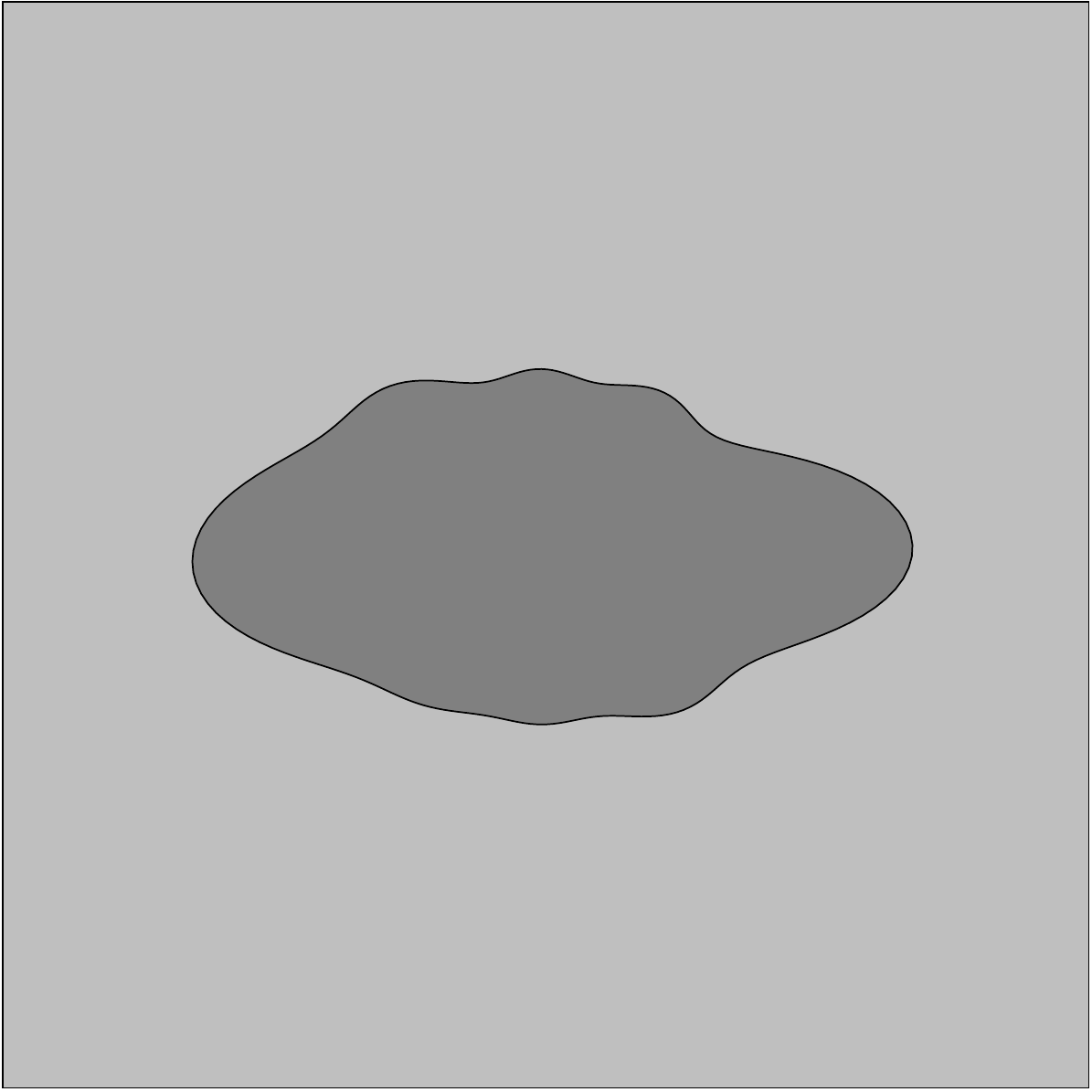}
\includegraphics[width=3cm,height=3cm]{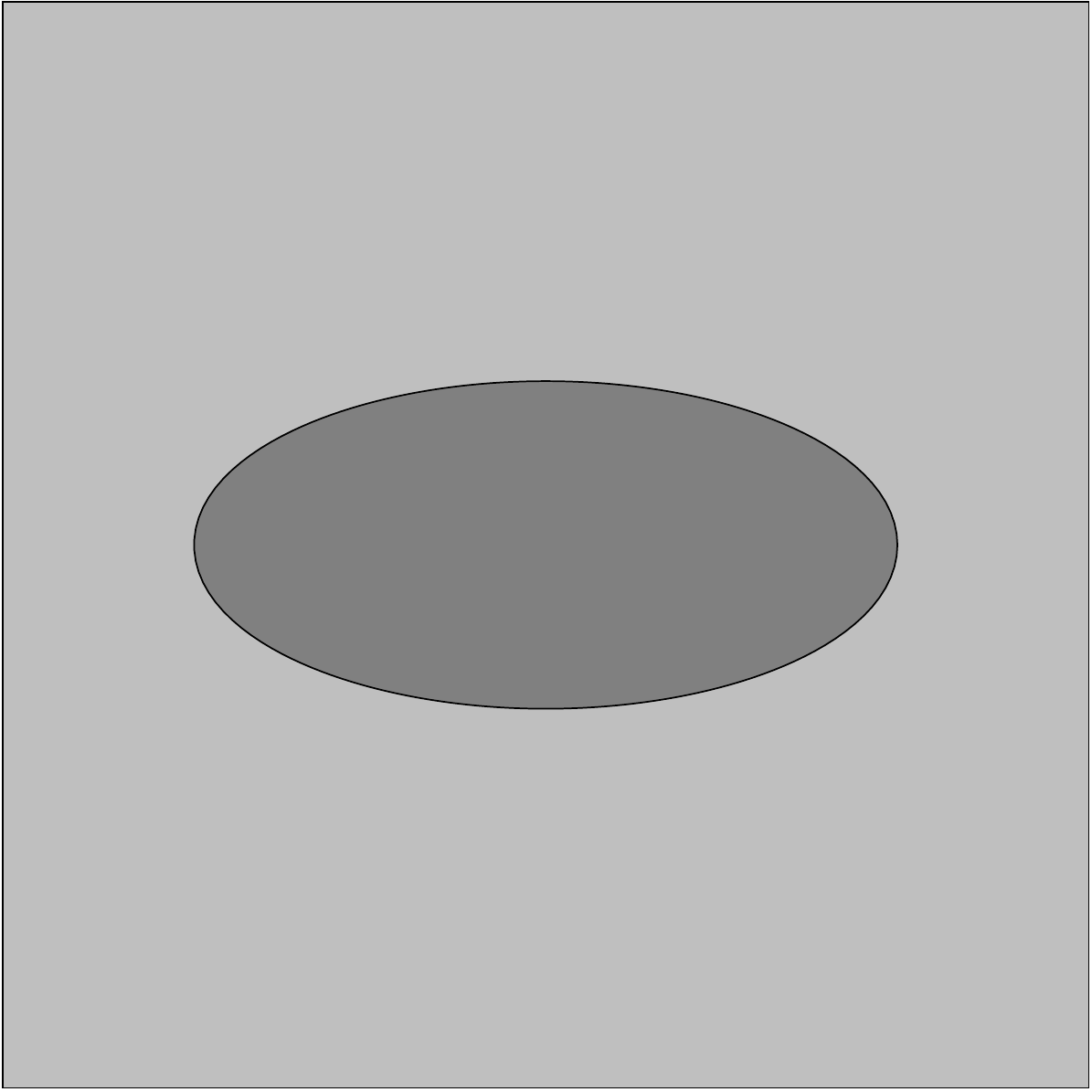}
\caption{\label{fig:shapes6}Different shapes which all generate 
the desired effective tensor $\boldsymbol{B}_5$ in case of 
a perforated domain.}
\end{center}
\end{figure}

We should finally have also a look at a perforated plate,
which has been considered Section~\ref{sec:scaffolds}.
We set $\sigma_1 = 1$ and choose
\[
\boldsymbol{B}_5 = \begin{bmatrix} 0.8 & 0 \\ 0 & 0.6 \end{bmatrix}.
\]
If we start the gradient method from the circle 
of radius $1/4$, we obtain the shape found in the
outermost left plot of Figure~\ref{fig:shapes6}. In 
case of starting with a randomly perturbed circle
of radius $1/4$, we get the two shapes seen in 
the middle plots in Figure~\ref{fig:shapes6}. Whereas, 
in the outermost right plot of Figure~\ref{fig:shapes6}, 
we plotted the an ellipse which also leads to the desired
effective tensor.

\section{Conclusion}
In this article, shape sensitivity analysis of the effective
tensor in case of composite materials and scaffold structures 
has been performed. In particular, we computed the shape gradient
of the least square matching of a desired material property.
In case of scaffold structures, we also provided the shape 
Hessian. This enables to apply the second order perturbation
method to quantify uncertainties in the geometry of the hole. 
Numerical tests based on a finite element implementation in
the two-dimensional setting have been performed for composite 
materials and for perforated domains. It has turned out that 
the computed optimal shapes depend on the initial guess, 
which means that the solution of the problem under 
consideration is in general not unique.
 

\end{document}